\documentclass[reqno,10pt, centertags,draft]{amsart}
\usepackage{amsmath,amsthm,amscd,amssymb,latexsym,upref}

\makeatletter
\def\theequation{\@arabic\c@equation}

\newcommand{\diag}{\operatorname{diag}}

\newcommand{\bbN}{{\mathbb{N}}}
\newcommand{\bbR}{{\mathbb{R}}}

\newcommand{\bbC}{{\mathbb{C}}}

\newcommand{\bbF}{{\mathbb{F}}} 
\newcommand{\bbH}{{\mathbb{H}}}

\newcommand{\cB}{{\mathcal B}}

\newcommand{\cF}{{\mathcal F}}
\newcommand{\cH}{{\mathcal H}}

\newcommand{\cK}{{\mathcal K}}
\newcommand{\cL}{{\mathcal L}}
\newcommand{\cM}{{\mathcal M}}
\newcommand{\cN}{{\mathcal N}}

\newcommand{\cS}{{\mathcal S}}

\newcommand{\no}{\nonumber}
\newcommand{\lb}{\label}
\newcommand{\f}{\frac}
\newcommand{\ul}{\underline}
\newcommand{\ol}{\overline}

\newcommand{\wti}{\widetilde  }

\newcommand{\oh}{o}

\newcommand{\loc}{\text{\rm{loc}}}
\newcommand{\tr}{\text{\rm{tr}}}
\newcommand{\spec}{\text{\rm{spec}}}

\newcommand{\dom}{\text{\rm{dom}}}

\newcommand{\supp}{\text{\rm{supp}}}

\newcommand{\bi}{\bibitem}

\newcommand{\hatt}{\widehat}
\renewcommand{\Re}{\text{\rm Re}}
\renewcommand{\Im}{\text{\rm Im}}

\renewcommand{\diag}{\text{\rm diag}}

\numberwithin{equation}{section}

\newtheorem{theorem}{Theorem}[section]
\newtheorem{lemma}[theorem]{Lemma}

\newtheorem{hypothesis}[theorem]{Hypothesis}
\newtheorem{remark}[theorem]{Remark}

\begin{document}

\title[Fredholm Determinants and Semi-Separable  
Kernels]{(Modified) Fredholm Determinants for Operators with
Matrix-Valued Semi-Separable Integral Kernels Revisited}
\author[F.\ Gesztesy and K.\ A.\ Makarov]{Fritz Gesztesy and Konstantin
A.\ Makarov}
\address{Department of Mathematics,
University of Missouri, Columbia, MO 65211, USA}
\email{fritz@math.missouri.edu}
\urladdr{http://www.math.missouri.edu/people/fgesztesy.html}
\address{Department of Mathematics, University of
Missouri, Columbia, MO 65211, USA}
\email{makarov@math.missouri.edu} 
\urladdr{http://www.math.missouri.edu/people/kmakarov.html}
\dedicatory{Dedicated with great pleasure to Eduard R.\ Tsekanovskii on
the occasion of his 65th birthday.}
\date{\today}
\thanks{{\it Integr. Equ. Oper. Theory} {\bf 47}, 457--497 (2003).}
\subjclass{Primary: 47B10, 47G10, Secondary: 34B27, 34L40.}
\keywords{Fredholm determinants, semi-separable kernels, Jost functions,
transmission coefficients, Floquet discriminants, Day's formula.}

\begin{abstract}
We revisit the computation of ($2$-modified) Fredholm determinants for
operators  with matrix-valued semi-separable integral kernels. The
latter occur, for instance, in the form of Green's functions associated
with closed ordinary differential operators on arbitrary intervals on
the real line. Our approach determines the ($2$-modified) Fredholm
determinants in terms of solutions of closely associated Volterra
integral equations, and as a result offers a natural way to compute
such determinants. 

We illustrate our approach by identifying classical objects such as the Jost
function for half-line Schr\"odinger operators and the inverse transmission
coefficient for Schr\"odinger operators on the real line as Fredholm
determinants, and rederiving the well-known expressions for them in due
course. We also apply our formalism to Floquet theory of Schr\"odinger
operators, and upon identifying the connection between the Floquet
discriminant and underlying Fredholm determinants, we derive new
representations of the Floquet discriminant. 

Finally, we rederive the explicit formula for the $2$-modified Fredholm
determinant corresponding to a convolution integral operator, whose
kernel is associated with a symbol given by a rational function, in a
straghtforward manner. This determinant formula represents a
Wiener--Hopf analog of Day's formula for the determinant associated
with finite Toeplitz matrices generated by the Laurent expansion of a
rational function.
\end{abstract}

\maketitle

\section{Introduction} \lb{s1}

We offer a self-contained and elementary approach to the computation of
Fredholm and $2$-modified Fredholm determinants associated with
$m\times m$ matrix-valued,  semi-separable integral kernels on arbitrary
intervals $(a,b)\subseteq\bbR$ of the type
\begin{equation}
K(x,x')=\begin{cases} f_1(x)g_1(x'), & a<x'< x< b, \\ 
f_2(x)g_2(x'), & a<x<x'<b, \end{cases}  \lb{1.1}
\end{equation}
associated with the Hilbert--Schmidt operator $K$ in $L^2((a,b);
dx)^{m}$, $m\in\bbN$, 
\begin{equation}
 (Kf)(x)=\int_a^b dx'\, K(x,x')f(x'), \quad f\in L^2((a,b); dx)^{m},  
\lb{1.2}
\end{equation}
assuming 
\begin{equation}  
f_j\in L^2((a,b); dx)^{m\times n_j}, \; 
g_j \in L^2((a,b); dx)^{n_j\times m}, \quad n_j\in\bbN, \; j=1,2. \lb{1.3}
\end{equation}
We emphasize that Green's matrices and resolvent operators associated with
closed ordinary differential operators on arbitrary intervals (finite or
infinite) on the real line are always of the form \eqref{1.1}--\eqref{1.3}
(cf.\ \cite[Sect.\ XIV.3]{GGK90}), as are certain classes of convolution
operators (cf.\ \cite[Sect.\ XIII.10]{GGK90}).

To describe the approach of this paper we briefly recall the principal
ideas of the approach to $m\times m$ matrix-valued semi-separable integral
kernels in the monographs by Gohberg, Goldberg, and Kaashoek \cite[Ch.\
IX]{GGK90} and Gohberg, Goldberg, and Krupnik \cite[Ch.\ XIII]{GGK00}. It
consists in decomposing $K$ in \eqref{1.2} into a Volterra operator $H_a$ and a
finite-rank operator $QR$
\begin{equation}
K=H_a + QR, \lb{1.4}  
\end{equation}
where
\begin{align}
(H_af)(x)&=\int_a^x dx'\, H(x,x')f(x'),  \quad f\in
L^2((a,b); dx)^{m}, \lb{1.5} \\
H(x,x')&=f_1(x)g_1(x')-f_2(x)g_2(x'), \quad a<x'<x<b \lb{1.6}
\end{align} 
and 
\begin{align}
&Q\colon \bbC^{n_2}\mapsto L^2((a,b); dx)^m, \quad (Q\ul u)(x)=f_2(x)\ul u,
\quad \ul u\in\bbC^{n_2}, \lb{1.7} \\
&R\colon L^2((a,b); dx)^m \mapsto \bbC^{n_2}, \quad (Rf)=\int_a^b
dx'\,g_2(x')f(x'), \quad f\in L^2((a,b); dx)^m. \lb{1.8} 
\end{align}
Moreover, introducing 
\begin{equation}
C(x)=(f_1(x) \;\; f_2(x)), \quad 
B(x)=(g_1(x) \;\; -g_2(x))^\top \lb{1.9} 
\end{equation}
and
the $n\times n$ matrix $A$ ($n=n_1+n_2$)
\begin{equation}
A(x)=\begin{pmatrix} g_1(x)f_1(x) & g_1(x)f_2(x) \\ 
-g_2(x)f_1(x) & -g_2(x)f_2(x) \end{pmatrix}, \lb{1.10}
\end{equation}
one considers a particular nonsingular solution $U(\cdot, \alpha)$ of
the following first-order system of differential equations
\begin{equation}
U'(x,\alpha )=\alpha A(x)U(x,\alpha ) \, \text{ for a.e.\ $x\in(a,b)$ and 
$\alpha\in\bbC$} \lb{1.11}
\end{equation}
and obtains 
\begin{align}
(I-\alpha H_a)^{-1}&=I+\alpha J_a(\alpha) \, \text{ for all 
$\alpha\in\bbC$,} \lb{1.12} \\ 
(J_a(\alpha) f)(x)&=\int_a^x dx'\,
J(x,x',\alpha )f(x'), \quad f\in L^2((a,b); dx)^m, \lb{1.13} \\  
J(x,x',\alpha )&=C(x)U(x,\alpha )U(x',\alpha )^{-1}B(x'), \quad a<x'<x<b. 
\lb{1.14}  
\end{align}
Next, observing
\begin{equation}
I-\alpha K =(I-\alpha H_a)[I-\alpha (I-\alpha H_a)^{-1}QR] \lb{1.15} 
\end{equation}
and assuming that $K$ is a trace class operator,
\begin{equation}
K\in\cB_1(L^2((a,b);dx)^m), \lb{1.15a}
\end{equation}
one computes,
\begin{align}
\det(I-\alpha K)&=\det(I-\alpha H_a)\det(I-\alpha (I-\alpha H_a)^{-1}QR) \no
\\  
&=\det(I-\alpha (I-\alpha H_a)^{-1}QR) \no \\
&={\det}_{\bbC^{n_2}}(I_{n_2}-\alpha R(I-\alpha H_a)^{-1}Q). \lb{1.16}  
\end{align}
In particular, the Fredholm determinant of $I-\alpha K$ is reduced to a
finite-dimensional determinant induced by the finite rank operator $QR$ in
\eqref{1.4}. Up to this point we followed the treatment in \cite[Ch.\
IX]{GGK90}). Now we will depart from the presentation in \cite[Ch.\
IX]{GGK90} and 
\cite[Ch.\ XIII]{GGK00} that focuses on a solution $U(\cdot, \alpha)$ of
\eqref{1.11} normalized by $U(a,\alpha)=I_n$. The latter normalization 
is in general not satisfied for Schr\"odinger operators on a half-line
or on the whole real line possessing eigenvalues as discussed in
Section \ref{s4}.  

To describe our contribution to this circle of ideas we now introduce the
Volterra integral equations
\begin{align}
\begin{split}
\hat f_1(x,\alpha )&=f_1(x)-\alpha \int_x^b dx'\, 
H(x,x')\hat f_1(x',\alpha ), \\ 
\hat f_2(x,\alpha )&=f_2(x)+\alpha \int_a^x dx'\, 
H(x,x')\hat f_2(x',\alpha ), \quad \alpha\in\bbC \lb{1.18} 
\end{split} 
\end{align}    
with solutions $\hat f_j(\cdot,\alpha ) \in L^2((a,b); dx)^{m\times n_j}$,  
$j=1,2$, and note that the first-order $n\times n$ system of differential
equations \eqref{1.11} then permits the explicit particular solution
\begin{align}
&U(x,\alpha )=\begin{pmatrix} I_{n_1}-\alpha \int_x^b dx'\, g_1(x')\hat
f_1(x',\alpha ) & \alpha\int_a^x dx'\, g_1(x')\hat f_2(x',\alpha ) \\ 
\alpha\int_x^b dx'\, g_2(x')\hat f_1(x',\alpha ) & I_{n_2}-\alpha \int_a^x
dx'\, g_2(x')\hat f_2(x',\alpha ) \end{pmatrix},  \no \\
& \hspace*{9.2cm} x\in (a,b). \lb{1.19}  
\end{align}
Given \eqref{1.19}, one can supplement \eqref{1.16} by 
\begin{align}
\det(I-\alpha K)&={\det}_{\bbC^{n_2}}(I_{n_2}-\alpha R(I-\alpha H_a)^{-1}Q)
\no \\
& ={\det}_{\bbC^{n_2}}\bigg(I_{n_2}-\alpha \int_a^b dx\,
g_2(x)\hat f_2(x,\alpha )\bigg) \no \\
& ={\det}_{\bbC^{n}}(U(b,\alpha)), \lb{1.20} 
\end{align}
our principal result. A similar set of results can of course be obtained
by introducing the corresponding Volterra operator $H_b$ in \eqref{2.5}.
Moreover, analogous results hold for $2$-modified Fredholm determinants
in the case where $K$ is only assumed to be a Hilbert--Schmidt
operator. 

Equations \eqref{1.16} and \eqref{1.20} summarize this approach based on
decomposing $K$ into a Volterra operator plus finite rank operator in
\eqref{1.4}, as advocated in \cite[Ch.\ IX]{GGK90} and \cite[Ch.\ XIII]{GGK00},
and our additional twist of relating this formalism to the underlying
Volterra integral equations \eqref{1.18} and the explicit solution \eqref{1.19}
of \eqref{1.11}. 

In Section \ref{s2} we set up the basic formalism leading up to the solution
$U$ in \eqref{1.19} of the first-order system of differential equations 
\eqref{1.11}. In Section \ref{s3} we derive the set of formulas \eqref{1.16},
\eqref{1.20}, if $K$ is a trace class operator, and their counterparts for
$2$-modified Fredholm determinants, assuming $K$ to be a
Hilbert--Schmidt operator only. Section \ref{s4} then treats four
particular applications: First we treat the case of half-line
Schr\"odinger operators in which we identify the Jost function as a
Fredholm determinant (a well-known, in fact, classical result due to
Jost and Pais \cite{JP51}). Next, we study the case of Schr\"odinger
operators on the real line in which we characterize the inverse of the
transmission coefficient as a Fredholm determinant (also a well-known
result, see, e.g., \cite[Appendix A]{Ne80}, \cite[Proposition
5.7]{Si79}). We also revisit this problem by replacing the second-order
Schr\"odinger equation by the equivalent first-order $2\times 2$ system
and determine the associated $2$-modified Fredholm determinant. The case
of periodic Schr\"odinger operators in which we derive a new
one-parameter family of representations of the Floquet discriminant and
relate it to underlying Fredholm determinants is discussed next.
Apparently, this is a new result. In our final Section \ref{s5}, we
rederive the explicit formula for the $2$-modified Fredholm determinant
corresponding to a convolution integral operator whose kernel is
associated with a symbol given by a rational function. The latter
represents a Wiener--Hopf analog of Day's formula \cite{Da75} for the
determinant of finite Toeplitz matrices generated by the Laurent
expansion of a rational function. The approach to ($2$-modified)
Fredholm determinants of semi-separable kernels advocated in this paper
permits a remarkably elementary derivation of this formula compared to
the current ones in the literature (cf.\ the references provided at the
end of Section \ref{s5}).

The effectiveness of the approach pursued in this paper is demonstrated 
by the ease of the computations involved and by the unifying character
it takes on when applied to differential and convolution-type operators
in several different settings.

\section{Hilbert--Schmidt operators with semi-separable integral kernels}
\lb{s2}

In this section we consider Hilbert-Schmidt operators with
matrix-valued semi-separable integral kernels following Gohberg,
Goldberg, and Kaashoek \cite[Ch.\ IX]{GGK90} and Gohberg, Goldberg, and
Krupnik \cite[Ch.\ XIII]{GGK00} (see also \cite{GK84}). To set up the basic
formalism we introduce the following hypothesis assumed throughout this
section. 
 
\begin{hypothesis} \lb{h2.1}  
Let $-\infty\leq a<b\leq \infty$ and $m,n_1,n_2\in\bbN$. Suppose that
$f_j$ are $m\times n_j$ matrices and $g_j$ are $n_j\times m$ matrices,
$j=1,2$, with $($Lebesgue$)$ measurable entries on $(a,b)$ such that
\begin{equation}  
f_j\in L^2((a,b); dx)^{m\times n_j}, \; 
g_j \in L^2((a,b); dx)^{n_j\times m}, \quad j=1,2. \lb{2.1}
\end{equation}
\end{hypothesis}

Given Hypothesis \ref{h2.1}, we introduce the Hilbert--Schmidt operator 
\begin{align}
\begin{split}
& K\in\cB_2(L^2((a,b); dx)^{m}), \lb{2.2} \\
& (Kf)(x)=\int_a^b dx'\, K(x,x')f(x'), \quad f\in L^2((a,b); dx)^{m}
\end{split}
\end{align}
in $L^2((a,b); dx)^{m}$ with $m\times m$ matrix-valued integral kernel
$K(\cdot,\cdot)$ defined by
\begin{equation}
K(x,x')=\begin{cases} f_1(x)g_1(x'), & a<x'< x< b, \\ 
f_2(x)g_2(x'), & a<x<x'<b. \end{cases}  \lb{2.3}
\end{equation}
One verifies that $K$ is a finite rank operator in $L^2((a,b); dx)^{m}$
if
$f_1=f_2$ and $g_1=g_2$ a.e. Conversely, any finite rank operator in
$L^2((a,b)); dx)^{m}$ is of the form \eqref{2.2}, \eqref{2.3}
with $f_1=f_2$ and $g_1=g_2$ (cf.\ \cite[p.\ 150]{GGK90}). 

Associated with $K$ we also introduce the Volterra operators $H_a$ and
$H_b$ in $L^2((a,b); dx)^{m}$ defined by 
\begin{align}
(H_af)(x)&=\int_a^x dx'\, H(x,x')f(x'), \lb{2.4} \\
(H_bf)(x)&=-\int_x^b dx'\, H(x,x')f(x'); \quad f\in
L^2((a,b); dx)^{m}, \lb{2.5} 
\end{align}
with $m\times m$ matrix-valued (triangular) integral kernel
\begin{equation}
H(x,x')=f_1(x)g_1(x')-f_2(x)g_2(x').  \lb{2.6}
\end{equation} 
Moreover, introducing the matrices\footnote{$M^\top$ denotes the
transpose of the matrix $M$.} 
\begin{align}
C(x)&=(f_1(x) \;\; f_2(x)), \lb{2.7} \\
B(x)&=(g_1(x) \;\; -g_2(x))^\top, \lb{2.8} 
\end{align}
one verifies
\begin{equation}
H(x,x')=C(x)B(x'), \, \text{ where } \begin{cases} 
a<x'<x<b & \text{for $H_a$,} \\ 
a<x<x'<b & \text {for $H_b$} \end{cases} \lb{2.9}
\end{equation}
and\footnote{$I_k$ denotes the identity matrix in $\bbC^k$, $k\in\bbN$.}
\begin{equation}
K(x,x')=\begin{cases} C(x)(I_n-P_0)B(x'), & a<x'<x<b, \\
-C(x)P_0B(x'), & a<x<x'<b \end{cases} \lb{2.9a}
\end{equation}
with 
\begin{equation}
 P_0=\begin{pmatrix} 0 & 0 \\ 0 & I_{n_2} \end{pmatrix}. \lb{2.9b} 
\end{equation}

Next, introducing the linear maps
\begin{align}
&Q\colon \bbC^{n_2}\mapsto L^2((a,b); dx)^m, \quad (Q\ul u)(x)=f_2(x)\ul u,
\quad \ul u\in\bbC^{n_2}, \lb{2.10} \\
&R\colon L^2((a,b); dx)^m \mapsto \bbC^{n_2}, \quad (Rf)=\int_a^b
dx'\,g_2(x')f(x'), \quad f\in L^2((a,b); dx)^m, \lb{2.11} \\
&S\colon \bbC^{n_1} \mapsto L^2((a,b); dx)^m, \quad (S \ul v)(x)=f_1(x)\ul v,
\quad \ul v\in\bbC^{n_1}, \lb{2.12} \\
&T\colon L^2((a,b); dx)^m \mapsto \bbC^{n_1}, \quad (Tf)=\int_a^b 
dx'\, g_1(x')f(x'), \quad f\in L^2((a,b); dx)^m, \lb{2.13} 
\end{align}
one easily verifies the following elementary yet significant result.
\begin{lemma} [\cite{GGK90}, Sect.\ IX.2; \cite{GGK00}, Sect.\ XIII.6]
\lb{l2.2} Assume Hypothesis \ref{h2.1}. Then
\begin{align}
K&=H_a + QR \lb{2.14} \\
 &=H_b + ST. \lb{2.15}
\end{align}
In particular, since $R$ and $T$ are of finite rank, so are $K-H_a$ and
$K-H_b$.
\end{lemma}

\begin{remark} \lb{r2.3}
The decompositions \eqref{2.14} and \eqref{2.15} of $K$ are significant
since they prove that $K$ is the sum of a Volterra and a finite rank
operator. As a consequence, the $($$2$-modified$)$ determinants
corresponding to $I-\alpha K$ can be reduced to determinants of
finite-dimensional matrices, as will be further discussed in Sections
\ref{s3} and \ref{s4}. 
\end{remark}

To describe the 
inverse\footnote{$I$ denotes the identity operator in $L^2((a,b); dx)^m$.} 
of $I-\alpha H_a$ and $I-\alpha H_b$,
$\alpha\in\bbC$, one introduces the $n\times n$ matrix $A$ ($n=n_1+n_2$)
\begin{align}
A(x)&=\begin{pmatrix} g_1(x)f_1(x) & g_1(x)f_2(x) \\ 
-g_2(x)f_1(x) & -g_2(x)f_2(x) \end{pmatrix} \lb{2.16} \\
& = B(x)C(x)\, \text{ for a.e.\ $x\in (a,b)$} \lb{2.17}
\end{align}
and considers a particular nonsingular solution $U=U(x,\alpha )$ of the
first-order $n\times n$  system of differential equations
\begin{equation}
U'(x,\alpha )=\alpha A(x)U(x,\alpha ) \, \text{ for a.e.\ $x\in(a,b)$ and 
$\alpha\in\bbC$.} \lb{2.18}
\end{equation}
Since $A\in L^1((a,b))^{n\times n}$, the general solution $V$ of
\eqref{2.18} is an $n\times n$ matrix with locally absolutely continuous
entries on $(a,b)$ of the form $V=UD$ for any constant  
$n\times n$ matrix $D$  (cf.\ \cite[Lemma IX.2.1]{GGK90})\footnote{If
$a>-\infty$, $V$ extends to an absolutely continuous $n\times n$ matrix on
all intervals of the type $[a,c)$, $c<b$. The analogous consideration
applies to the endpoint $b$ if $b<\infty$.}.

\begin{theorem} [\cite{GGK90}, Sect.\ IX.2; \cite{GGK00}, Sects.\ XIII.5,
XIII.6] \lb{t2.4} ${}$ \\
Assume Hypothesis~\ref{h2.1} and let $U(\cdot,\alpha)$
denote a nonsingular solution of \eqref{2.18}. Then, \\
$(i)$ $I-\alpha H_a$ and $I-\alpha H_b$ are invertible for all
$\alpha\in\bbC$ and 
\begin{align}
(I-\alpha H_a)^{-1}&=I+\alpha J_a(\alpha), \lb{2.19} \\
(I-\alpha H_b)^{-1}&=I+\alpha J_b(\alpha), \lb{2.20} \\
(J_a(\alpha) f)(x)&=\int_a^x dx'\, J(x,x',\alpha )f(x'), \lb{2.21} \\ 
(J_b(\alpha) f)(x)&=-\int_x^b dx'\, J(x,x',\alpha )f(x'); \quad f\in
L^2((a,b); dx)^m, \lb{2.22} \\  
J(x,x',\alpha )&=C(x)U(x,\alpha )U(x',\alpha )^{-1}B(x'), \, \text{ where }
\begin{cases}  a<x'<x<b & \text{for $J_a$,} \\ 
a<x<x'<b & \text {for $J_b$.} \end{cases}  \lb{2.23}  
\end{align}
$(ii)$ Let $\alpha\in\bbC$. Then $I-\alpha K$ is invertible if and only if
the
$n_2\times n_2$ matrix $I_{n_2}-\alpha R(I-\alpha H_a)^{-1}Q$ is. Similarly,
$I-\alpha K$ is invertible if and only if the $n_1\times n_1$ matrix
$I_{n_1}-\alpha T(I-\alpha H_b)^{-1}S$ is. In particular,
\begin{align}
(I-\alpha K)^{-1}&=(I-\alpha H_a)^{-1}+\alpha (I-\alpha H_a)^{-1}QR(I-\alpha K)^{-1} \lb{2.24} \\
&=(I-\alpha H_a)^{-1} \no \\
& \quad +\alpha (I-\alpha H_a)^{-1}Q[I_{n_2}-\alpha R(I-\alpha H_a)^{-1}Q]^{-1}R(I-\alpha H_a)^{-1}
\lb{2.25} \\
&=(I-\alpha H_b)^{-1}+\alpha (I-\alpha H_b)^{-1}ST(I-\alpha K)^{-1} \lb{2.26} \\
&=(I-\alpha H_b)^{-1} \no \\
& \quad +\alpha (I-\alpha H_b)^{-1}S[I_{n_1}-\alpha T(I-\alpha H_b)^{-1}S]^{-1}T(I-\alpha H_b)^{-1}.
\lb{2.27}
\end{align}
Moreover,
\begin{align}
(I-\alpha K)^{-1}&=I+\alpha L(\alpha), \lb{2.28} \\
(L(\alpha) f)(x)&=\int_a^b dx'\, L(x,x',\alpha )f(x'), \lb{2.29} \\  
L(x,x',\alpha )&=\begin{cases} C(x)U(x,\alpha )(I-P(\alpha))
U(x',\alpha)^{-1}B(x'), & a<x'<x<b, \\  
-C(x)U(x,\alpha )P(\alpha)U(x',\alpha )^{-1}B(x'), & a<x<x'<b, \end{cases} 
\lb{2.30}  
\end{align}
where $P(\alpha)$ satisfies
\begin{equation}
P_0U(b,\alpha )(I-P(\alpha))=(I-P_0)U(a,\alpha )P(\alpha), \quad 
P_0=\begin{pmatrix} 0 & 0 \\ 0 & I_{n_2} \end{pmatrix}. \lb{2.31}
\end{equation}
\end{theorem}

\begin{remark} \lb{r2.5}  
$(i)$ The results \eqref{2.19}--\eqref{2.23} and
\eqref{2.28}--\eqref{2.31} are easily verified by computing
$(I-\alpha H_a)(I+\alpha J_a)$ and $(I+\alpha J_a)(I-\alpha H_a)$, etc., using an integration
by parts. Relations \eqref{2.24}--\eqref{2.27} are clear from
\eqref{2.14} and \eqref{2.15}, a standard resolvent identity, and the
fact that $K-H_a$ and $K-H_b$ factor into $QR$ and $ST$, respectively. \\
$(ii)$ The discussion in \cite[Sect.\ IX.2]{GGK90}, \cite[Sects.\
XIII.5, XIII.6]{GGK00} starts from the particular normalization 
\begin{equation}
U(a,\alpha )=I_n \lb{2.32}
\end{equation}
of a solution $U$ satisfying \eqref{2.18}. In this case the
explicit solution for $P(\alpha)$ in \eqref{2.31} is given by
\begin{equation}
P(\alpha)=\begin{pmatrix} 0 & 0 \\ U_{2,2}(b,\alpha )^{-1}
U_{2,1}(b,\alpha ) & I_{n_2} \end{pmatrix}. \lb{2.33}
\end{equation}
However, for concrete applications to differential operators to be discussed
in Section \ref{s4}, the normalization \eqref{2.32} is not necessarily
possible.
\end{remark}

Rather than solving the basic first-order system of differential equations
$U'=\alpha AU$ in \eqref{2.18} with the fixed initial condition 
$U(a,\alpha)=I_n$ in \eqref{2.32}, we now derive an explicit particular solution
of \eqref{2.18} in terms of closely associated solutions of Volterra integral
equations involving the integral kernel $H(\cdot,\cdot)$ in \eqref{2.6}. This
approach is most naturally suited for the applications to Jost functions,
transmission coefficients, and Floquet discriminants we discuss in
Section \ref{s4} and to the class of Wiener--Hopf operators we study in
Section \ref{s5}. 

Still assuming Hypothesis \ref{h2.1}, we now introduce the Volterra
integral equations
\begin{align}
\hat f_1(x,\alpha )&=f_1(x)-\alpha \int_x^b dx'\, H(x,x')\hat f_1(x',\alpha ), \lb{2.35} \\ 
\hat f_2(x,\alpha )&=f_2(x)+\alpha \int_a^x dx'\, H(x,x')\hat f_2(x',\alpha ); 
\quad \alpha\in\bbC, \lb{2.36}  
\end{align}  
with solutions $\hat f_j(\cdot,\alpha ) \in L^2((a,b); dx)^{m\times n_j}$,
$j=1,2$.
\begin{lemma} \lb{l2.6}
Assume Hypothesis \ref{h2.1} and let $\alpha\in\bbC$. \\
$(i)$ The first-order $n\times n$ system of differential equations
$U'=\alpha AU$ a.e.\ on $(a,b)$ in \eqref{2.18} permits the explicit particular
solution
\begin{align}
&U(x,\alpha )=\begin{pmatrix} I_{n_1}-\alpha \int_x^b dx'\, g_1(x')\hat
f_1(x',\alpha ) & \alpha\int_a^x dx'\, g_1(x')\hat f_2(x',\alpha ) \\ 
\alpha\int_x^b dx'\, g_2(x')\hat f_1(x',\alpha ) & I_{n_2}-\alpha \int_a^x
dx'\, g_2(x')\hat f_2(x',\alpha ) \end{pmatrix},  \no \\
& \hspace*{9.2cm} x\in (a,b). \lb{2.37}  
\end{align}
As long as\footnote{${\det}_{\bbC^{k}}(M)$ and $\tr_{\bbC^{k}}(M)$ denote
the determinant and trace of a $k\times k$ matrix $M$ with complex-valued
entries, respectively.} 
\begin{align}
& {\det}_{\bbC^{n_1}}\bigg(I_{n_1}-\alpha \int_a^b dx\, g_1(x)\hat
f_1(x,\alpha )\bigg) \neq 0, \lb{2.38} \\
\intertext{or equivalently,} 
&{\det}_{\bbC^{n_2}}\bigg(I_{n_2}-\alpha \int_a^b dx\, 
g_2(x)\hat f_2(x,\alpha )\bigg) \neq 0, \lb{2.39} 
\end{align}
$U$ is nonsingular for all $x\in (a,b)$ and the general
solution $V$ of \eqref{2.18} is then of the form $V=UD$ for any constant 
$n\times n$ matrix $D$. \\
$(ii)$ Choosing \eqref{2.37} as the particular solution $U$ in
\eqref{2.28}--\eqref{2.31}, $P(\alpha)$ in \eqref{2.31} simplifies to 
\begin{equation}
P(\alpha)=P_0=\begin{pmatrix} 0 & 0
\\ 0 & I_{n_2} \end{pmatrix}. \lb{2.40}
\end{equation}
\end{lemma}
\begin{proof}
Differentiating the right-hand side of \eqref{2.37} with respect to $x$ and 
using the Volterra integral equations \eqref{2.35}, \eqref{2.36} readily
proves that $U$ satisfies $U'=\alpha AU$ a.e.\ on $(a,b)$. 

By Liouville's formula (cf., e.g., \cite[Theorem IV.1.2]{Ha82}) one infers 
\begin{equation}
{\det}_{\bbC^n}(U(x,\alpha))={\det}_{\bbC^n}(U(x_0,\alpha))
\exp\bigg(\alpha \int_{x_0}^x dx'\, \tr_{\bbC^n}(A(x'))\bigg), \quad 
x,x_0 \in (a,b). \lb{2.41} 
\end{equation}
Since $\tr_{\bbC^n}(A)\in L^1((a,b);dx)$ by \eqref{2.1},
\begin{equation}
\lim_{x\downarrow a} {\det}_{\bbC^n}(U(x,\alpha)) \, \text{ and } \,  
\lim_{x\uparrow b} {\det}_{\bbC^n}(U(x,\alpha)) \, \text{ exist.} 
\lb{2.43} 
\end{equation}
Hence, if \eqref{2.38} holds, $U(x,\alpha)$ is nonsingular for $x$ in a
neighborhood $(a,c)$, $a<c$,  of $a$, and similarly, if \eqref{2.39} holds,
$U(x,\alpha)$ is nonsingular for $x$ in a neighborhood $(c,b)$, $c<b$, of $b$.
In either case, \eqref{2.41} then proves that $U(x,\alpha)$ is nonsingular for
all $x\in (a,b)$. 

Finally, since $U_{2,1}(b,\alpha )=0$, \eqref{2.40} follows from \eqref{2.33}.
\end{proof}

\begin{remark} \lb{r2.7} 
In concrete applications $($e.g., to Schr\"odinger operators on a
half-line or on the whole real axis as discussed in Section \ref{s4}$)$,
it may happen that ${\det}_{\bbC^n}(U(x,\alpha))$ vanishes for certain
values of intrinsic parameters $($such as the energy parameter$)$.
Hence, a normalization of the type $U(a,\alpha)=I_n$  is impossible in 
the case of such parameter values and the normalization  of $U$ is best
left open as illustrated in Section
\ref{s4}. One also  observes that in general our explicit particular 
solution $U$ in \eqref{2.37} satisfies $U(a,\alpha)\neq I_n$,
$U(b,\alpha)\neq I_n$. 
\end{remark}

\begin{remark} \lb{r2.8} 
In applications to Schr\"odinger and Dirac-type systems, $A$ is
typically of the form 
\begin{equation}
A(x)=e^{-Mx}\wti A(x) e^{Mx}, \quad x\in(a,b) \lb{2.44}
\end{equation}
where $M$ is an $x$-independent $n\times n$ matrix $($in general
depending on a spectral parameter$)$ and $\wti A$ has a simple
asymptotic behavior such that for some $x_0\in(a,b)$ 
\begin{equation}
\int_a^{x_0} w_a(x) dx\, |\wti A(x)-\wti A_-| + 
\int_{x_0}^b w_b(x) dx\, |\wti A(x)-\wti A_+| < \infty
\lb{2.45}
\end{equation}
for constant $n\times n$ matrices $\wti A_{\pm}$ and appropriate
weight functions $w_a\geq 0$, $w_b\geq 0$. Introducing
$W(x,\alpha)=e^{Mx} U(x,\alpha)$, equation \eqref{2.18} reduces to
\begin{equation}
W'(x,\alpha)=[M+\alpha \wti A(x)]W(x,\alpha), \quad x\in(a,b) \lb{2.46}
\end{equation}
with
\begin{equation}
{\det}_{\bbC^n} (W(x,\alpha))={\det}_{\bbC^n} (U(x,\alpha)) 
e^{-\tr_{\bbC^n}(M)x}, \quad x\in(a,b). \lb{2.47}
\end{equation}
The system \eqref{2.46} then leads to operators $H_a$, $H_b$, and $K$. 
We will briefly illustrate this in connection with Schr\"odinger
operators on the line in Remark \ref{r4.5b}.
\end{remark}

\section{(Modified) Fredholm determinants for operators with
semi-separable integral kernels} \lb{s3}

In the first part of this section we suppose that $K$ is a trace class
operator and consider the Fredholm determinant of $I-K$. In the second
part we consider $2$-modified Fredholm determinants in the case where
$K$ is a Hilbert--Schmidt operator.

In the context of trace class operators we assume the following
hypothesis. 

\begin{hypothesis} \lb{h3.1}  
In addition to Hypothesis \ref{h2.1}, we suppose that $K$ is a trace
class operator, $K\in\cB_1(L^2((a,b); dx)^m)$. 
\end{hypothesis}

The following results can be found in Gohberg, Goldberg, and
Kaashoek \cite[Theorem 3.2]{GGK90} and in Gohberg, Goldberg, and Krupnik
\cite[Sects.\ XIII.5, XIII.6]{GGK00} under the additional assumptions that
$a,b$ are finite and $U$ satisfies the normalization $U(a)=I_n$ (cf.\
\eqref{2.18}, \eqref{2.32}). Here we present the general case where
$(a,b)\subseteq\bbR$ is an arbitrary interval on the real line and $U$ is not
normalized but given by the particular solution \eqref{2.37}.

In the course of the proof we use some of the standard properties of
determinants, such as,
\begin{align}
& \det((I_\cH-A)(I_\cH-B))=\det(I_\cH-A)\det(I_\cH-B), \quad A, B
\in\cB_1(\cH), \lb{3.0} \\   
& \det(I_{\cH_1}-AB)=\det(I_\cH-BA) \;\, \text{ for all
$A\in\cB_1(\cH_1,\cH)$, $B\in\cB(\cH,\cH_1)$} \lb{3.0a} \\
& \hspace*{5cm} \text{ such that $AB\in\cB_1(\cH_1)$, $BA\in \cB_1(\cH)$,} \no
\intertext{and}
&\det(I_\cH-A)={\det}_{\bbC^k}(I_k-D_k) \, \text{ for } \, A=\begin{pmatrix} 
0 & C \\ 0 & D_k \end{pmatrix}, \;\, \cH=\cK\dotplus \bbC^{k}, \lb{3.0b} 
\intertext{since}
&I_\cH -A=\begin{pmatrix} I_\cK & -C \\ 0 & I_k-D_k \end{pmatrix} =
\begin{pmatrix} I_\cK & 0 \\ 0 & I_k-D_k \end{pmatrix} 
\begin{pmatrix} I_\cK & -C \\ 0 & I_k \end{pmatrix}. \lb{3.0c} 
\end{align}
Here $\cH$ and $\cH_1$ are complex separable Hilbert spaces,
$\cB(\cH)$ denotes the set of bounded linear operators on $\cH$, $\cB_p(\cH)$,
$p\geq 1$, denote the usual trace ideals of $\cB(\cH)$, and $I_\cH$ denotes
the identity operator in $\cH$. Moreover, $\det_p(I_\cH-A)$, $A\in\cB_p(\cH)$,
denotes the ($p$-modified) Fredholm determinant of $I_\cH-A$ with
$\det_1(I_\cH-A)=\det(I_\cH-A)$, $A\in\cB_1(\cH)$, the standard Fredholm
determinant of a trace class operator, and $\tr(A)$, $A\in\cB_1(\cH)$, the
trace of a trace class operator. Finally, $\dotplus$ in \eqref{3.0b} denotes a
direct but not necessary orthogonal direct decomposition of $\cH$ into $\cK$
and the $k$-dimensional subspace $\bbC^k$. (We refer, e.g., to \cite{GGK96},
\cite[Sect.\ IV.1]{GK69}, \cite[Ch.\ 17]{RS78}, \cite{Si77}, 
\cite[Ch.\ 3]{Si79} for these facts). 
 
\begin{theorem} \lb{t3.2} Suppose Hypothesis~\ref{h3.1} and let
$\alpha\in\bbC$. Then, 
\begin{align} 
&\tr(H_a)=\tr(H_b)=0, \quad \det(I-\alpha H_a)=\det(I-\alpha H_b)=1, 
\lb{3.1} \\ 
&\tr(K)=\int_a^b dx\,\tr_{\bbC^{n_1}}(g_1(x)f_1(x))
=\int_a^b dx\,\tr_{\bbC^{m}}(f_1(x)g_1(x)) \lb{3.2} \\
& \hspace*{.86cm} =\int_a^b dx\,\tr_{\bbC^{n_2}}(g_2(x)f_2(x))
=\int_a^b dx\,\tr_{\bbC^{m}}(f_2(x)g_2(x)). \lb{3.3} 
\end{align}
Assume in addition that $U$ is given by \eqref{2.37}. Then,
\begin{align}
&\det(I-\alpha K)={\det}_{\bbC^{n_1}}(I_{n_1}-\alpha T(I
-\alpha H_b)^{-1}S) \lb{3.4}
\\ 
& \hspace*{1.95cm} ={\det}_{\bbC^{n_1}}\bigg(I_{n_1}-\alpha \int_a^b dx\,
g_1(x)\hat f_1(x,\alpha )\bigg) \lb{3.5} \\
& \hspace*{1.95cm} ={\det}_{\bbC^{n}}(U(a,\alpha)) \lb{3.5a} \\
& \hspace*{1.95cm} ={\det}_{\bbC^{n_2}}(I_{n_2}-\alpha R(I
-\alpha H_a)^{-1}Q)
\lb{3.6} \\
& \hspace*{1.95cm} ={\det}_{\bbC^{n_2}}\bigg(I_{n_2}-\alpha \int_a^b dx\,
g_2(x)\hat f_2(x,\alpha )\bigg) \lb{3.7} \\
& \hspace*{1.95cm} ={\det}_{\bbC^{n}}(U(b,\alpha)). \lb{3.8} 
\end{align}
\end{theorem}
\begin{proof}
We briefly sketch the argument following \cite[Theorem 3.2]{GGK90} since
we use a different solution $U$ of $U'=\alpha AU$. Relations \eqref{3.1}
are clear from Lidskii's theorem (cf., e.g., \cite[Theorem
VII.6.1]{GGK90}, \cite[Sect.\ III.8, Sect.\ IV.1]{GK69}, \cite[Theorem
3.7]{Si79}). Thus, 
\begin{equation}
\tr(K)=\tr(QR)=\tr(RQ)=\tr(ST)=\tr(TS) \lb{3.9}
\end{equation}
then proves \eqref{3.2} and \eqref{3.3}. Next, one observes
\begin{align}
I-\alpha K &=(I-\alpha H_a)[I-\alpha (I-\alpha H_a)^{-1}QR] \lb{3.10} \\ 
&=(I-\alpha H_b)[I-\alpha (I-H_b)^{-1}ST] \lb{3.11}
\end{align}
and hence,
\begin{align}
\det(I-\alpha K)&=\det(I-\alpha H_a)\det(I-\alpha (I-\alpha H_a)^{-1}QR)
\no \\ 
&=\det(I-\alpha (I-\alpha H_a)^{-1}QR) \no \\
&=\det(I-\alpha R(I-\alpha H_a)^{-1}Q) \no \\
&={\det}_{\bbC^{n_2}}(I_{n_2}-\alpha R(I-\alpha H_a)^{-1}Q) \lb{3.12} \\
&={\det}_{\bbC^{n}}(U(b,\alpha)). \lb{3.13} 
\end{align}
Similarly,
\begin{align} 
\det(I-\alpha K)&=\det(I-\alpha H_b)\det(I-\alpha (I-\alpha H_b)^{-1}ST)
\no \\ 
&=\det(I-\alpha (I-\alpha H_b)^{-1}ST) \no \\
&=\det(I-\alpha T(I-\alpha H_b)^{-1}S) \no \\
&={\det}_{\bbC^{n_1}}(I_{n_1}-\alpha T(I-\alpha H_b)^{-1}S) \lb{3.14} \\
&={\det}_{\bbC^{n}}(U(a,\alpha)). \lb{3.15} 
\end{align}
Relations \eqref{3.13} and \eqref{3.15} follow directly from 
taking the limit $x\uparrow b$ and $x\downarrow a$ in \eqref{2.38}. This 
proves \eqref{3.4}--\eqref{3.8}.
\end{proof}

Equality of \eqref{3.13} and \eqref{3.15} also follows directly from 
\eqref{2.41} and 
\begin{align}
\int_a^b dx' \tr_{\bbC^n}(A(x'))&=\int_a^b dx'
[\tr_{\bbC^{n_1}}(g_1(x')f_1(x'))-\tr_{\bbC^{n_2}}(g_2(x')f_2(x'))] \lb{3.17a}
\\ 
&=\tr(H_a)=\tr(H_b)=0. \lb{3.18}
\end{align}

Finally, we treat the case of $2$-modified Fredholm determinants in the
case  where $K$ is only assumed to lie in the Hilbert-Schmidt class. In
addition to
\eqref{3.0}--\eqref{3.0b} we will use the following standard facts for
$2$-modified Fredholm determinants $\det_2(I-A)$, $A\in\cB_2(\cH)$
(cf.\, e,g., 
\cite{GGK97}, \cite[Ch.\ XIII]{GGK00}, \cite[Sect.\ IV.2]{GK69},
\cite{Si77}, \cite[Ch.\ 3]{Si79}), 
\begin{align}
& {\det}_2 (I-A)= \det((I-A)\exp(A)), \quad A\in\cB_2 (\cH), \lb{3.21} \\
& {\det}_2((I-A)(I-B))={\det}_2(I-A){\det}_2(I-B)e^{-\tr(AB)}, \quad A, B
\in\cB_2(\cH), \lb{3.22} \\
& {\det}_2(I-A)=\det(I-A)e^{\tr(A)}, \quad A\in\cB_1(\cH). \lb{3.23}
\end{align}

\begin{theorem} \lb{t3.3} Suppose Hypothesis~\ref{h2.1} and let
$\alpha\in\bbC$. Then,
\begin{equation}
{\det}_2 (I-\alpha H_a)={\det}_2 (I-\alpha H_b)=1. \lb{3.24} 
\end{equation}
Assume in addition that $U$ is given by \eqref{2.37}. Then,
\begin{align}
&{\det}_2 (I-\alpha K)={\det}_{\bbC^{n_1}}(I_{n_1}-\alpha T(I
-\alpha H_b)^{-1}S) \exp(\alpha \, {\tr}_{\bbC^m}(ST)) \lb{3.25} \\   
& \hspace*{1.95cm} ={\det}_{\bbC^{n_1}}\bigg(I_{n_1}-\alpha \int_a^b dx\,
g_1(x)\hat f_1(x,\alpha )\bigg) \no \\
& \hspace*{2.4cm} \times \exp\bigg(\alpha \int_a^b dx\,
{\tr}_{\bbC^m}(f_1(x)g_1(x))\bigg) \lb{3.26} \\
& \hspace*{1.95cm} ={\det}_{\bbC^{n}}(U(a,\alpha)) 
\exp\bigg(\alpha \int_a^b
dx\, {\tr}_{\bbC^m}(f_1(x)g_1(x))\bigg) \lb{3.27} \\
& \hspace*{1.95cm} ={\det}_{\bbC^{n_2}}(I_{n_2}-\alpha R(I
-\alpha H_a)^{-1}Q) \exp(\alpha \, {\tr}_{\bbC^m}(QR)) \lb{3.28} \\
& \hspace*{1.95cm} ={\det}_{\bbC^{n_2}}\bigg(I_{n_2}-\alpha \int_a^b dx\,
g_2(x)\hat f_2(x,\alpha )\bigg) \no \\
& \hspace*{2.4cm} \times \exp\bigg(\alpha \int_a^b dx\,
{\tr}_{\bbC^m}(f_2(x)g_2(x))\bigg) \lb{3.29} \\
& \hspace*{1.95cm} ={\det}_{\bbC^{n}}(U(b,\alpha)) \exp\bigg(\alpha \int_a^b
dx\, {\tr}_{\bbC^m}(f_2(x)g_2(x))\bigg). \lb{3.30} 
\end{align}
\end{theorem}
\begin{proof}
Relations \eqref{3.24} follow since the Volterra operators $H_a, H_b$ have
no nonzero eigenvalues. Next, again using \eqref{3.10} and \eqref{3.11}, 
one computes,
\begin{align}
{\det}_2 (I-\alpha K) 
& ={\det}_2 (I-\alpha H_a){\det}_2 (I-\alpha (I-\alpha H_a)^{-1}QR) \no
\\  
& \qquad \times \exp(-\tr(\alpha^2 H_a(I-\alpha H_a)^{-1}QR)) \no \\
& =\det(I-\alpha (I-\alpha H_a)^{-1}QR) 
\exp(\alpha \, \tr((I-\alpha H_a)^{-1}QR))\no \\
& \quad \times \exp(-\tr(\alpha^2 H_a(I-\alpha H_a)^{-1}QR)) \no \\
& ={\det}_{\bbC^{n_2}}(I_{n_2}-\alpha R(I-\alpha H_a)^{-1}Q) 
\exp(\alpha \, \tr(QR)) \lb{3.32} \\
& ={\det}_{\bbC^{n}}(U(b,\alpha))\exp\bigg(\alpha \int_a^b dx\,
{\tr}_{\bbC^m}(f_1(x)g_1(x))\bigg). \lb{3.33} 
\end{align}
Similarly,
\begin{align}
{\det}_2 (I-\alpha K) 
& ={\det}_2 (I-\alpha H_b){\det}_2 (I-\alpha (I-\alpha H_b)^{-1}ST) \no \\  
& \quad \times \exp(-\tr(\alpha^2 H_b(I-\alpha H_b)^{-1}ST)) \no \\
& =\det(I-\alpha (I-\alpha H_b)^{-1}ST)  
\exp(\alpha \, \tr((I-\alpha H_b)^{-1}ST))\no \\
& \quad \times \exp(-\tr(\alpha^2 H_b(I-\alpha H_b)^{-1}ST)) \no \\
& ={\det}_{\bbC^{n_1}}(I_{n_1}-\alpha T(I-\alpha H_b)^{-1}S) 
\exp(\alpha \, \tr(ST))\lb{3.34} \\ 
& ={\det}_{\bbC^{n}}(U(a,\alpha))\exp\bigg(\alpha
\int_a^b dx\, {\tr}_{\bbC^m}(f_2(x)g_2(x))\bigg). \lb{3.35} 
\end{align}
\end{proof}

Equality of \eqref{3.33} and \eqref{3.35} also follows directly from 
\eqref{2.41} and \eqref{3.17a}.

\section{Some applications to Jost functions, transmission coefficients,
and Floquet discriminants of Schr\"odinger operators} \lb{s4}

In this section we illustrate the results of Section \ref{s3} in 
three particular cases: The case of Jost functions for half-line Schr\"odinger
operators, the transmission coefficient for Schr\"odinger operators on the
real line,  and the case of Floquet discriminants associated with
Schr\"odinger  operators on a compact interval. The case of a the
second-order Schr\"odinger operator on the line is also transformed
into a first-order $2 \times 2$ system and its associated $2$-modified
Fredholm deteminant is identified with that of the Schr\"odinger
operator on $\bbR$.  For simplicity we will limit ourselves to scalar
coefficients although the results for half-line Schr\"odinger operators
and those on the full real line immediately extend to the matrix-valued
situation.  

\medskip

We start with the case of half-line Schr\"odinger operators:\\
{\bf The case $\boldsymbol{(a,b)=(0,\infty)}$:} Assuming
\begin{equation}
V\in L^1((0,\infty); dx), 
\lb{4.1}
\end{equation}
(we note that $V$ is not necessarily assumed to be real-valued) we
introduce the closed Dirichlet-type operators in $L^2((0,\infty); dx)$
defined by 
\begin{align}
&H^{(0)}_+f=-f'', \no \\
&f\in \dom\big(H^{(0)}_+\big)=\{g\in L^2((0,\infty); dx) \,|\, g,g' \in
AC_{\loc}([0,R]) 
\text{ for all $R>0$}, \lb{4.3} \\
& \hspace*{2.95cm} f(0_+)=0, \, f''\in L^2((0,\infty); dx)\}, \no \\
&H_+f=-f''+Vf, \no \\
&f\in\dom(H_+)=\{g\in L^2((0,\infty); dx) \,|\, g,g' \in AC_{\loc}([0,R]) 
\text{ for all $R>0$}, \lb{4.4} \\
& \hspace*{2.7cm} f(0_+)=0, \, (-f''+Vf) \in L^2((0,\infty); dx)\}. \no
\end{align}
We note that $H^{(0)}_+$ is self-adjoint and that $H_+$ is self-adjoint
if and only if $V$ is real-valued. 

Next we introduce the regular solution $\phi(z,\cdot)$ and Jost solution 
$f(z,\cdot)$ of $-\psi''(z)+V\psi(z)=z\psi(z)$, $z\in\bbC\backslash\{0\}$, by
\begin{align}
\phi(z,x)&=z^{-1/2}\sin(z^{1/2}x)+\int_0^x dx' g^{(0)}_+(z,x,x')
V(x')\phi(z,x'), \lb{4.5} \\
f(z,x)&=e^{iz^{1/2}x}-\int_x^\infty dx' g^{(0)}_+(z,x,x')V(x')f(z,x'), \lb{4.6}
\\
&\hspace*{2.72cm} \Im(z^{1/2})\geq 0, \; z\neq 0, \; x\geq 0, \no 
\end{align}
where 
\begin{align}
g^{(0)}_+(z,x,x')&=z^{-1/2}\sin(z^{1/2}(x-x')).  \lb{4.7}
\end{align}
We also introduce the Green's function of $H^{(0)}_+$,
\begin{equation}
G^{(0)}_+(z,x,x')=\big(H^{(0)}_+-z\big)^{-1}(x,x')=\begin{cases}  
z^{-1/2}\sin(z^{1/2}x)e^{iz^{1/2}x'}, & x\leq x', \\
z^{-1/2}\sin(z^{1/2}x')e^{iz^{1/2}x}, & x \geq x'. \lb{4.8}
\end{cases}
\end{equation}
The Jost function $\cF$ associated with the pair $\big(H_+,H^{(0)}_+\big)$ is
given by
\begin{align}
\cF(z)&=W(f(z),\phi(z))=f(z,0) \lb{4.9} \\
&=1+z^{-1/2}\int_0^\infty dx\, \sin(z^{1/2}x)V(x)f(z,x) \lb{4.10} \\
&=1+\int_0^\infty dx\, e^{iz^{1/2}x}V(x)\phi(z,x); \quad 
 \Im(z^{1/2})\geq 0, \; z\neq 0, \lb{4.11} 
\end{align}
where 
\begin{equation}
W(f,g)(x)=f(x)g'(x)-f'(x)g(x), \quad x \geq 0, \lb{4.12}
\end{equation}
denotes the Wronskian of $f$ and $g$. Introducing the factorization
\begin{equation}
V(x)=u(x)v(x), \quad u(x)=|V(x)|^{1/2}\exp(i\arg(V(x))), \;
v(x)=|V(x)|^{1/2}, \lb{4.13}
\end{equation}
one verifies\footnote{$\ol T$ denotes the operator closure of $T$ and
$\spec(\cdot)$ abbreviates the spectrum of a linear operator.}
\begin{align}
(H_+-z)^{-1}&=\big(H^{(0)}_+-z\big)^{-1} \no \\
& \quad -\big(H^{(0)}_+-z\big)^{-1}v
\Big[I+\ol{u\big(H^{(0)}_+-z\big)^{-1}v}\Big]^{-1}u\big(H^{(0)}_+-z\big)^{-1}, 
\lb{4.14}  \\ 
& \hspace*{6.4cm} z\in\bbC\backslash\spec(H_+). \no 
\end{align}
To establish the connection with the notation used in Sections \ref{s2}
and \ref{s3}, we introduce the operator $K(z)$ in $L^2((0,\infty); dx)$
(cf.\ \eqref{2.3}) by
\begin{equation}
K(z)=-\ol{u\big(H^{(0)}_+-z\big)^{-1}v}, \quad z\in\bbC\backslash
\spec\big(H^{(0)}_+\big)
\lb{4.15}
\end{equation}
with integral kernel
\begin{equation}
K(z,x,x')=-u(x)G^{(0)}_+(z,x,x')v(x'), \quad \Im(z^{1/2})\geq 0, 
\; x,x' \geq 0, \lb{4.16} 
\end{equation}
and the Volterra operators $H_0 (z)$, $H_\infty (z)$ (cf.\ \eqref{2.4},
\eqref{2.5}) with integral kernel
\begin{equation}
H(z,x,x')=u(x)g^{(0)}_+(z,x,x')v(x').  \lb{4.17} 
\end{equation}
Moreover, we introduce for a.e.\ $x > 0$, 
\begin{align}
\begin{split}
f_1(z,x)&=-u(x) e^{iz^{1/2}x}, \hspace*{1.82cm} 
g_1(z,x)=v(x)z^{-1/2}\sin(z^{1/2}x),  \\
f_2(z,x)&=-u(x)z^{-1/2}\sin(z^{1/2}x), \quad g_2(z,x)=v(x) e^{iz^{1/2}x}. 
\lb{4.18}
\end{split}
\end{align}
Assuming temporarily that 
\begin{equation}
\supp(V) \text{ is compact} \lb{4.20}
\end{equation}
in addition to hypothesis \eqref{4.1}, introducing $\hat f_j(z,x)$, $j=1,2$, by
\begin{align}
\hat f_1(z,x)&=f_1(z,x)-\int_x^\infty dx'\, H(z,x,x')\hat f_1(z,x'), \lb{4.21}
\\ 
\hat f_2(z,x)&=f_2(z,x)+\int_0^x dx'\, H(z,x,x')\hat f_2(z,x'), \lb{4.22} \\
&\hspace*{1.85cm}  \Im(z^{1/2})\geq 0, \; z\neq 0, \; x \geq 0, \no
\end{align}  
yields solutions $\hat f_j(z,\cdot) \in L^2((0,\infty); dx)$, $j=1,2$. By
comparison with \eqref{4.5}, \eqref{4.6}, one then identifies
\begin{align}
\hat f_1(z,x)&=-u(x) f(z,x), \lb{4.23} \\
\hat f_2(z,x)&=-u(x) \phi(z,x). \lb{4.24} 
\end{align}
We note that the temporary compact support assumption \eqref{4.20} on $V$ has
only been introduced to guarantee that $f_2(z,\cdot), \hat f_2(z,\cdot) \in
L^2((0,\infty); dx)$. This extra hypothesis will soon be removed. 

We start with a well-known result.

\begin{theorem} [Cf., e.g., \cite{RS79}, Theorem XI.20] \lb{t4.2}
Suppose $f,g \in L^q(\bbR; dx)$ for some $2\leq q<\infty$. Denote by $f(X)$
the maximally defined multiplication operator by $f$ in $L^2(\bbR; dx)$ and by
$g(P)$ the maximal multiplication operator by $g$ in Fourier 
space\footnote{That is, $P=-id/dx$ with domain $\dom(P)=H^{2,1}(\bbR)$, the
usual Sobolev space.} $L^2(\bbR; dp)$. Then\footnote{$\cB_q(\cH)$, $q\geq 1$
denote the usual trace ideals, cf.\ \cite{GK69}, \cite{Si79}. }  
$f(X)g(P)\in \cB_q(L^2(\bbR; dx))$ and
\begin{equation}
\|f(X)g(P)\|_{\cB_q(L^2(\bbR; dx))} \leq (2\pi)^{-1/q} 
\|f\|_{L^q(\bbR; dx)}\|g\|_{L^q(\bbR; dx)}. \lb{4.24a}
\end{equation}
\end{theorem}

We will use Theorem \ref{t4.2}, to sketch a proof of the following known
result:

\begin{theorem} \lb{t4.1}
Suppose $V\in L^1((0,\infty); dx)$ and let $z\in\bbC$ with
$\Im(z^{1/2})>0$. Then 
\begin{equation}
K(z)\in \cB_1(L^2((0,\infty); dx)). \lb{4.24A}
\end{equation}
\end{theorem}
\begin{proof}
For $z<0$ this is discussed in the proof of \cite[Theorem XI.31]{RS79}.
For completeness we briefly sketch the principal arguments of a proof of
Theorem
\ref{t4.1}. One possible approach consists of reducing Theorem \ref{t4.1} to
Theorem \ref{t4.2} in the special case $q=2$ by embedding the half-line problem
on $(0,\infty)$ into a problem on $\bbR$ as follows. One introduces the
decomposition
\begin{equation}
L^2(\bbR; dx)=L^2((0,\infty); dx) \oplus L^2((-\infty,0); dx), \lb{4.24b} 
\end{equation} 
and extends $u,v,V$ to $(-\infty,0)$ by putting $u,v,V$ equal to zero on
$(-\infty,0)$, introducing
\begin{equation}
\tilde u(x)=\begin{cases} u(x), & x>0, \\ 0, & x<0, \end{cases} \quad 
\tilde v(x)=\begin{cases} v(x), & x>0, \\ 0, & x<0, \end{cases} \quad
\widetilde V(x)=\begin{cases} V(x), & x>0, \\ 0, & x<0. \end{cases} \lb{4.24ba}
\end{equation}
Moreover, one considers the Dirichlet Laplace operator $H^{(0)}_D$ in 
$L^2(\bbR;dx)$ by
\begin{align}
H^{(0)}_Df&=-f'', \no \\
\dom\big(H^{(0)}_D\big)&=\{g\in L^2(\bbR; dx) \,|\, g,g' \in
AC_{\loc}([0,R])\cap AC_{\loc}([-R,0]) \text{ for all $R>0$}, \no \\
& \hspace*{.65cm} f(0_\pm)=0, \, f''\in L^2(\bbR; dx)\} \lb{4.24c} 
\end{align}
and introduces
\begin{equation}
\widetilde K(z)=-\ol{\tilde u \big(H^{(0)}_D -z\big)^{-1} \tilde v} 
=K(z)\oplus 0, \quad \Im(z^{1/2})>0. \lb{4.24ca}
\end{equation}
By Krein's formula, the resolvents of the Dirichlet Laplace operator
$H^{(0)}_D$ and that of the ordinary Laplacian $H^{(0)}=P^2=-d^2/dx^2$ on
$H^{2,2}(\bbR)$ differ precisely by a rank one operator. Explicitly, one
obtains
\begin{align} 
G^{(0)}_D(z,x,x')&=G^{(0)}(z,x,x')-G^{(0)}(z,x,0)G^{(0)}(z,0,0)^{-1}
G^{(0)}(z,0,x') \no \\
&=G^{(0)}(z,x,x')-\f{i}{2z^{1/2}}\exp(iz^{1/2}|x|)\exp(iz^{1/2}|x'|),
\lb{4.24cb} \\
& \hspace*{4.25cm} \Im(z^{1/2})>0, \; x,x'\in\bbR, \no
\end{align}
where we abbreviated the Green's functions of $H^{(0)}_D$ and
$H^{(0)}=-d^2/dx^2$ by
\begin{align}
G^{(0)}_D(z,x,x')&=\big(H^{(0)}_D -z\big)^{-1}(x,x'), \lb{4.24cC} \\
G^{(0)}(z,x,x')&=\big(H^{(0)}-z\big)^{-1}(x,x')=
\f{i}{2z^{1/2}}\exp(iz^{1/2}|x-x'|). \lb{4.24cc}
\end{align}
Thus,
\begin{equation}
\widetilde K(z)=-\ol{\tilde u\big(H^{(0)}-z\big)^{-1}\tilde v} 
-\f{i}{2z^{1/2}} \big(\tilde v \, 
\ol{\exp(iz^{1/2}|\cdot|)}\,,\,\cdot\,\big)\tilde u \exp(iz^{1/2}|\cdot|).
\lb{4.24cd}
\end{equation}
By Theorem \ref{t4.2} for $q=2$ one infers that 
\begin{equation}
\big[\tilde u\big(H^{(0)}-z\big)^{-1/2}\big] \in 
\cB_2(L^2(\bbR;dx)), \quad \Im(z^{1/2})>0 \lb{4.24ce}
\end{equation}
and hence,
\begin{equation}
\big[\tilde u\big(H^{(0)}-z\big)^{-1/2}\big] 
\ol{\big[\big(H^{(0)}-z\big)^{-1/2}\tilde v\big]} \in 
\cB_1(L^2(\bbR;dx)), \quad \Im(z^{1/2})>0. \lb{4.24cf}
\end{equation}
Since the second term on the right-hand side of \eqref{4.24cd} is a rank one
operator one concludes
\begin{equation}
\widetilde K(z) \in \cB_1(L^2(\bbR;dx)), \quad \Im(z^{1/2})>0 \lb{4.24cg}
\end{equation}
and hence \eqref{4.24A} using \eqref{4.24ca}.
\end{proof}

An application of Lemma \ref{l2.6} and Theorem \ref{t3.2} then yields the
following well-known result identifying the Fredholm determinant of $I-K(z)$
and the Jost function $\cF(z)$.

\begin{theorem} \lb{t4.3}
Suppose $V\in L^1((0,\infty); dx)$ and let $z\in\bbC$ with
$\Im(z^{1/2})>0$. Then 
\begin{equation}
\det(I-K(z))=\cF(z). \lb{4.25}
\end{equation}
\end{theorem}
\begin{proof}
Assuming temporarily that $\supp(V)$ is compact (cf.\ \eqref{4.20}), Lemma
\ref{l2.6}  applies and one obtains from \eqref{2.37} and
\eqref{4.18}--\eqref{4.24} that
\begin{align}
U(z,x)&=\begin{pmatrix} 1- \int_x^\infty dx'\, g_1(z,x')\hat
f_1(z,x') & \int_0^x dx'\, g_1(z,x')\hat f_2(z,x') \\[1mm] 
\int_x^\infty dx'\, g_2(z,x')\hat f_1(z,x') & 1-\int_0^x
dx'\, g_2(z,x')\hat f_2(z,x') \end{pmatrix},  \no \\
&=\begin{pmatrix} 1+ \int_x^\infty dx'\, \f{\sin(z^{1/2}x')}{z^{1/2}}V(x')
f(z,x') & -\int_0^x dx'\, \f{\sin(z^{1/2}x')}{z^{1/2}}V(x') \phi(z,x')
\\[1mm]  
-\int_x^\infty dx'\, e^{iz^{1/2}x'} V(x') f(z,x') & 1+\int_0^x
dx'\, e^{iz^{1/2}x'} V(x') \phi(z,x') \end{pmatrix}, \no \\
& \hspace*{9cm} x>0. \lb{4.27} 
\end{align}
Relations \eqref{3.5} and \eqref{3.7} of Theorem \ref{t3.2} with
$m=n_1=n_2=1$, $n=2$, then immediately yield
\begin{align}
\det(I-K(z))&=1+z^{-1/2}\int_0^\infty dx\, \sin(z^{1/2}x)V(x)f(z,x) \no
\\ &= 1+\int_0^\infty dx\, e^{iz^{1/2}x}V(x)\phi(z,x) \no \\
&=\cF(z) \lb{4.28}
\end{align}
and hence \eqref{4.25} is proved under the additional hypothesis \eqref{4.20}.
Removing the compact support hypothesis on $V$ now follows by a standard
argument. For completeness we sketch this argument next. Multiplying
$u,v,V$ by a smooth cutoff function $\chi_\varepsilon$ of compact support of
the type
\begin{equation}
0\leq \chi\leq 1, \quad \chi(x)=\begin{cases} 1, & x\in [0,1], \\ 
0, & |x|\geq 2, \end{cases} \quad \chi_\varepsilon (x)=\chi(\varepsilon x), \; 
\varepsilon >0, \lb{4.29}
\end{equation}
denoting the results by $u_\varepsilon=u\chi_\varepsilon$,
$v_\varepsilon=v\chi_\varepsilon$, $V_\varepsilon=V\chi_\varepsilon$, one
introduces in analogy to \eqref{4.24ba},
\begin{equation}
\tilde u_\varepsilon (x)=\begin{cases} u_\varepsilon (x), & x>0, \\ 0, & x<0,
\end{cases} \quad 
\tilde v_\varepsilon (x)=\begin{cases} v_\varepsilon (x), & x>0, \\ 0, & x<0,
\end{cases} \quad
\tilde V_\varepsilon (x)=\begin{cases} V_\varepsilon (x), & x>0, \\ 0, & x<0,
\end{cases} \lb{4.29a}
\end{equation}
and similarly, in analogy to \eqref{4.15} and \eqref{4.24ca},
\begin{align}
K_\varepsilon (z)&=-\ol{u_\varepsilon  \big(H^{(0)}_+-z\big)^{-1}
v_\varepsilon}, \quad \Im(z^{1/2})>0, \lb{4.29b} \\
\widetilde K_\varepsilon (z)&=-\ol{\tilde u_\varepsilon  \big(H^{(0)}_D
-z\big)^{-1} \tilde v_\varepsilon}  =K_\varepsilon (z)\oplus 0, 
\quad \Im(z^{1/2})>0. \lb{4.29c}
\end{align}
One then estimates,
\begin{align}
& \big\|\widetilde K(z)-\widetilde K_\varepsilon (z)
\big\|_{\cB_1(L^2(\bbR;dx))} \no \\  
& \quad \leq \Big\|-\ol{\tilde u\big(H^{(0)}-z\big)^{-1}\tilde v}
+\ol{\tilde u_\varepsilon \big(H^{(0)}-z\big)^{-1}\tilde v_\varepsilon} 
\Big\|_{\cB_1(L^2(\bbR;dx))} \no \\
& \qquad\, + \f{1}{2|z|^{1/2}} 
\Big\|\big(\tilde v \, \ol{\exp(iz^{1/2}|\cdot|)}\, , \, \cdot \,\big) 
\tilde u \exp(iz^{1/2}|\cdot|)  \no \\
& \hspace*{2.35cm} \, -
\big(\tilde v_\varepsilon \, \ol{\exp(iz^{1/2}|\cdot|)}\,,\,\cdot\,\big) 
\tilde u_\varepsilon  \exp(iz^{1/2}|\cdot|)\Big\|_{\cB_1(L^2(\bbR;dx))} \no \\
& \quad \leq \Big\|-\ol{\tilde u\big(H^{(0)}-z\big)^{-1}\tilde v}
+\ol{\tilde u_\varepsilon \big(H^{(0)}-z\big)^{-1}\tilde v} \no \\
& \hspace*{1cm} -\ol{\tilde u_\varepsilon \big(H^{(0)}-z\big)^{-1}\tilde v}  
+\ol{\tilde u_\varepsilon \big(H^{(0)}-z\big)^{-1}\tilde v_\varepsilon} 
\Big\|_{\cB_1(L^2(\bbR;dx))} \no \\
& \qquad\, + \f{1}{2|z|^{1/2}} 
\Big\|\big(\tilde v \, \ol{\exp(iz^{1/2}|\cdot|)}\, , \, \cdot \,\big) 
\tilde u \exp(iz^{1/2}|\cdot|)  \no \\
& \hspace*{2.35cm} \,
-\big(\tilde v \, \ol{\exp(iz^{1/2}|\cdot|)}\,,\,\cdot\,\big) 
\tilde u_\varepsilon \exp(iz^{1/2}|\cdot|) \no \\
& \hspace*{2.35cm} \,
+\big(\tilde v \, \ol{\exp(iz^{1/2}|\cdot|)}\,,\,\cdot\,\big) 
\tilde u_\varepsilon \exp(iz^{1/2}|\cdot|) \no \\
& \hspace*{2.35cm} \,
-\big(\tilde v_\varepsilon \, \ol{\exp(iz^{1/2}|\cdot|)}\,,\,\cdot\,\big) 
\tilde u_\varepsilon  \exp(iz^{1/2}|\cdot|)\Big\|_{\cB_1(L^2(\bbR;dx))} 
\no 
\end{align}
\begin{align}
& \hspace*{1.1cm} \leq \widetilde C(z) 
\big[\|\tilde u-\tilde u_\varepsilon\|_{L^2(\bbR;dx)}
+ \|\tilde v-\tilde v_\varepsilon\|_{L^2(\bbR;dx)}\big]  
= C(z) \|\tilde v-\tilde v_\varepsilon\|_{L^2(\bbR;dx)} \no \\
& \hspace*{1.1cm} \leq C(z) \|v-v_\varepsilon\|_{L^2((0,\infty);dx)},
\lb{4.29d}
\end{align}
where $C(z)=2\widetilde C(z)>0$ is an appropriate constant. Thus, applying
\eqref{4.24ca} and \eqref{4.29c}, one finally concludes 
\begin{equation}
\lim_{\varepsilon\downarrow 0}\big\|K(z)-K_\varepsilon (z)
\big\|_{\cB_1(L^2((0,\infty); dx))}=0. \lb{4.31} 
\end{equation}
Since $V_\varepsilon$ has compact support, \eqref{4.28} applies to
$V_\varepsilon$ and one obtains,
\begin{equation}
\det(I-K_\varepsilon (z))=\cF_\varepsilon (z), \lb{4.33}
\end{equation}
where, in obvious notation, we add the subscript $\varepsilon$ to all
quantities associated with $V_\varepsilon$ resulting in $\phi_\varepsilon$,
$f_\varepsilon$, $\cF_\varepsilon$, $f_{\varepsilon, j}$, 
$\hat f_{\varepsilon, j}$, $j=1,2$, etc. By \eqref{4.31}, the left-hand side of
\eqref{4.33} converges to $\det(I-K(z))$ as $\varepsilon\downarrow 0$. Since
\begin{equation}
\lim_{\varepsilon\downarrow 0}\|V_\varepsilon -V\|_{L^1((0,\infty); dx)}=0,
\lb{4.34}
\end{equation}
the Jost function $\cF_\varepsilon$ is well-known to converge to $\cF$
pointwise as $\varepsilon\downarrow 0$ (cf.\ \cite{BG83}). Indeed, fixing $z$
and iterating the  Volterra integral equation \eqref{4.6} for $f_\varepsilon$
shows that
$|z^{-1/2}\sin(z^{1/2}x)f_\varepsilon (z,x)|$ is uniformly bounded with respect
to $(x,\varepsilon)$ and hence the continuity of $\cF_\varepsilon (z)$ with
respect to $\varepsilon$ follows from \eqref{4.34} and the analog of
\eqref{4.10} for $V_\varepsilon$,
\begin{equation}
\cF_\varepsilon (z)=1+ z^{-1/2}\int_0^\infty dx\, \sin(z^{1/2}x)
V_\varepsilon (x) f_\varepsilon (z,x), \lb{4.35}
\end{equation}
applying the dominated convergence theorem. Hence, 
\eqref{4.33} yields \eqref{4.25} in the limit $\varepsilon\downarrow 0$.  
\end{proof} 

\begin{remark} \lb{r4.3a}  
$(i)$ The result \eqref{4.28} explicitly shows that ${\det}_{\bbC^n}(U(z,0))$
vanishes for each eigenvalue $z$ $($one then necessarily has $z<0$$)$ of the
Schr\"odinger operator $H$. Hence, a normalization of the type
$U(z,0)=I_n$ is clearly impossible in such a case. \\
$(ii)$ The right-hand side $\cF$ of \eqref{4.25} $($and hence the
Fredholm determinant on the left-hand side$)$ admits a continuous extension to
the positive real line. Imposing the additional exponential falloff of
the potential of the type $V\in L^1((0,\infty); \exp(ax)dx)$ for some $a>0$,
then $\cF$ and hence the Fredholm determinant on the left-hand side of 
\eqref{4.25} permit an analytic continuation through the essential spectrum of
$H_+$ into a strip of width $a/2$ $($w.r.t. the variable $z^{1/2}$$)$. This is
of particular relevance in the study of resonances of $H_+$ $($cf.\
\cite{Si00}$)$. 
\end{remark}

The result \eqref{4.25} is well-known, we refer, for instance, to \cite{JP51},
\cite{Ne61}, \cite{Ne72}, \cite[p.\ 344--345]{Ne02}, \cite{Si00}.
(Strictly speaking, these authors additionally assume $V$ to be
real-valued, but this is not essential in this context.) The current
derivation  presented appears to be by far the simplest available in
the literature as it only involves the elementary manipulations leading
to
\eqref{3.4}--\eqref{3.8}, followed by a standard approximation argument to
remove the compact support hypothesis on $V$.

Since one is dealing with the Dirichlet Laplacian on $(0,\infty)$ 
in the half-line context, Theorem \ref{t4.1} extends to a larger potential
class characterized by
\begin{equation}
\int_0^R dx\, x |V(x)| + \int_R^\infty dx\, |V(x)| < \infty \lb{4.36}
\end{equation}
for some fixed $R>0$. We omit the corresponding details but refer to
\cite[Theorem XI.31]{RS79},  which contains the necessary basic facts to make
the transition from hypothesis
\eqref{4.1} to \eqref{4.36}. 

\medskip

Next we turn to Schr\"odinger operators on the real line: \\
{\bf The case $\boldsymbol{(a,b)=\bbR}$:}  Assuming
\begin{equation}
V\in L^1(\bbR; dx), 
\lb{4.39}
\end{equation}
we introduce the closed operators in $L^2(\bbR; dx)$ defined by 
\begin{align}
&H^{(0)}f=-f'', \quad f\in \dom\big(H^{(0)}\big)=H^{2,2}(\bbR), \lb{4.40} \\
&Hf=-f''+Vf, \lb{4.41} \\
&f\in\dom(H)=\{g\in L^2(\bbR; dx) \,|\, g,g' \in AC_{\loc}(\bbR); \, (-f''+Vf)
\in L^2(\bbR); dx)\}.  \no
\end{align}
Again, $H^{(0)}$ is self-adjoint. Moreover, $H$ is self-adjoint if and
only if $V$ is real-valued.

Next we introduce the Jost solutions $f_\pm (z,\cdot)$ of
$-\psi''(z)+V\psi(z)=z\psi(z)$, $z\in\bbC\backslash\{0\}$, by 
\begin{align}
f_\pm (z,x)&=e^{\pm iz^{1/2}x}-\int_x^{\pm\infty} dx'
g^{(0)}(z,x,x')V(x')f_\pm (z,x'), \lb{4.42} \\
&\hspace*{3.34cm} \Im(z^{1/2})\geq 0, \; z\neq 0, \; x\in\bbR, \no 
\end{align}
where $g^{(0)}(z,x,x')$ is still given by \eqref{4.7}. We also introduce the
Green's function of $H^{(0)}$,
\begin{equation}
G^{(0)}(z,x,x')=\big(H^{(0)}-z\big)^{-1}(x,x')=
\f{i}{2z^{1/2}}e^{iz^{1/2}|x-x'|}, \quad \Im(z^{1/2})>0, \; x,x'\in\bbR.
\lb{4.43}
\end{equation}
The Jost function $\cF$ associated with the pair $\big(H,H^{(0)}\big)$ is given
by
\begin{align}
\cF(z)&=\f{W(f_-(z),f_+(z))}{2iz^{1/2}}  \lb{4.44} \\
&=1-\f{1}{2iz^{1/2}}\int_\bbR dx\, e^{\mp iz^{1/2}x}V(x)f_\pm (z,x), 
\quad \Im(z^{1/2})\geq 0, \; z\neq 0, \lb{4.45}
\end{align}
where $W(\cdot,\cdot)$ denotes the Wronskian defined in \eqref{4.12}. We note
that if $H^{(0)}$ and $H$ are self-adjoint, then
\begin{equation}
T(\lambda)=\lim_{\varepsilon\downarrow 0} \cF(\lambda+i\varepsilon)^{-1}, 
\quad \lambda>0,
\end{equation}
denotes the transmission coefficient corresponding to the pair
$\big(H,H^{(0)}\big)$. Introducing again the factorization \eqref{4.13} of
$V=uv$, one verifies as in \eqref{4.14} that 
\begin{align}
(H-z)^{-1}&=\big(H^{(0)}-z\big)^{-1} \no \\
& \quad -\big(H^{(0)}-z\big)^{-1}v
\Big[I+\ol{u\big(H^{(0)}-z\big)^{-1}v}\Big]^{-1}u\big(H^{(0)}-z\big)^{-1}, 
\lb{4.46} \\ 
& \hspace*{6.6cm} z\in\bbC\backslash\spec(H). \no
\end{align}
To make contact with the notation used in Sections \ref{s2} and \ref{s3}, we
introduce the operator $K(z)$ in $L^2(\bbR; dx)$ (cf.\ \eqref{2.3},
\eqref{4.15}) by
\begin{equation}
K(z)=-\ol{u\big(H^{(0)}-z\big)^{-1}v}, \quad
z\in\bbC\backslash\spec\big(H^{(0)}\big) \lb{4.47}
\end{equation}
with integral kernel
\begin{equation}
K(z,x,x')=-u(x)G^{(0)}(z,x,x')v(x'), \quad \Im(z^{1/2})\geq 0, \, z\neq 0, 
\; x,x' \in\bbR,  \lb{4.48} 
\end{equation}
and the Volterra operators $H_{-\infty} (z)$, $H_\infty (z)$ (cf.\ \eqref{2.4},
\eqref{2.5}) with integral kernel
\begin{equation}
H(z,x,x')=u(x)g^{(0)}(z,x,x')v(x').  \lb{4.49} 
\end{equation}
Moreover, we introduce for a.e.\ $x\in\bbR$, 
\begin{align}
\begin{split}
f_1(z,x)&=-u(x) e^{iz^{1/2}x}, \hspace*{.68cm} 
g_1(z,x)=(i/2)z^{-1/2}v(x)e^{-iz^{1/2}x},  \\
f_2(z,x)&=-u(x)e^{-iz^{1/2}x}, \quad \; 
g_2(z,x)=(i/2)z^{-1/2}v(x)e^{iz^{1/2}x}. \lb{4.50}
\end{split}
\end{align}
Assuming temporarily that 
\begin{equation}
\supp(V) \text{ is compact} \lb{4.52}
\end{equation}
in addition to hypothesis \eqref{4.39}, introducing $\hat f_j(z,x)$, $j=1,2$,
by
\begin{align}
\hat f_1(z,x)&=f_1(z,x)-\int_x^\infty dx'\, H(z,x,x')\hat f_1(z,x'), \lb{4.53}
\\ 
\hat f_2(z,x)&=f_2(z,x)+\int_{-\infty}^x dx'\, H(z,x,x')\hat f_2(z,x'),
\lb{4.54} \\ 
&\hspace*{2cm}  \Im(z^{1/2})\geq 0, \; z\neq 0, \; x \in\bbR, \no
\end{align}  
yields solutions $\hat f_j(z,\cdot) \in L^2(\bbR; dx)$, $j=1,2$. By
comparison with \eqref{4.42}, one then identifies
\begin{align}
\hat f_1(z,x)&=-u(x) f_+ (z,x), \lb{4.55} \\
\hat f_2(z,x)&=-u(x) f_- (z,x). \lb{4.56} 
\end{align}
We note that the temporary compact support assumption \eqref{4.20} on $V$ has
only been introduced to guarantee that $f_j(z,\cdot), \hat f_j(z,\cdot) \in
L^2(\bbR; dx)$, $j=1,2$. This extra hypothesis will soon be removed. 

We also recall the well-known result.

\begin{theorem} \lb{t4.4}
Suppose $V\in L^1(\bbR; dx)$ and let $z\in\bbC$ with $\Im(z^{1/2})>0$.
Then 
\begin{equation}
K(z)\in \cB_1(L^2(\bbR; dx)). \lb{4.57}
\end{equation}
\end{theorem}

This is an immediate consequence of Theorem \ref{t4.2} with $q=2$. 

An application of Lemma \ref{l2.6} and Theorem \ref{t3.2} then again yields the
following well-known result identifying the Fredholm determinant of $I-K(z)$
and the Jost function $\cF(z)$ (inverse transmission coefficient).

\begin{theorem} \lb{t4.5}
Suppose $V\in L^1(\bbR; dx)$ and let $z\in\bbC$ with $\Im(z^{1/2})>0$.
Then 
\begin{equation}
\det(I-K(z))=\cF(z). \lb{4.58}
\end{equation} 
\end{theorem}
\begin{proof}
Assuming temporarily that $\supp(V)$ is compact (cf.\ \eqref{4.20}), Lemma
\ref{l2.6}  applies and one infers from \eqref{2.37} and
\eqref{4.50}--\eqref{4.56} that
\begin{equation}
U(z,x)=\begin{pmatrix} 1- \int_x^\infty dx'\, g_1(z,x')\hat
f_1(z,x') & \int_{-\infty}^x dx'\, g_1(z,x')\hat f_2(z,x') \\ 
\int_x^\infty dx'\, g_2(z,x')\hat f_1(z,x') & 1-\int_{-\infty}^x
dx'\, g_2(z,x')\hat f_2(z,x') \end{pmatrix}, \quad x\in\bbR, \lb{4.59}
\end{equation}
becomes
\begin{align}
U_{1,1}(z,x)&=1+ \f{i}{2z^{1/2}}\int_x^\infty dx'\,
e^{-iz^{1/2}x'}V(x') f_+ (z,x'), \lb{4.60} \\
U_{1,2}(z,x)&= -\f{i}{2z^{1/2}}\int_{-\infty}^x dx'\,
e^{-iz^{1/2}x'}V(x') f_- (z,x'), \lb{4.61} \\  
U_{2,1}(z,x)&=-\f{i}{2z^{1/2}}\int_x^\infty dx'\, e^{iz^{1/2}x'} V(x') f_+
(z,x'), \lb{4.62} \\ 
U_{2,2}(z,x)&= 1+\f{i}{2z^{1/2}}\int_{-\infty}^x dx'\, e^{iz^{1/2}x'} V(x')f_-
(z,x'). \lb{4.63} 
\end{align}
Relations \eqref{3.5} and \eqref{3.7} of Theorem \ref{t3.2} with
$m=n_1=n_2=1$, $n=2$, then immediately yield
\begin{align}
\det(I-K(z))&=1-\f{1}{2iz^{1/2}}\int_\bbR dx\, e^{\mp iz^{1/2}x}V(x)f_\pm (z,x)
\no \\  
&=\cF(z) \lb{4.64}
\end{align}
and hence \eqref{4.58} is proved under the additional hypothesis \eqref{4.52}.
Removing the compact support hypothesis on $V$ now follows line by line the 
approximation argument discussed in the proof of Theorem \ref{t4.3}.  
\end{proof} 

Remark \ref{r4.3a} applies again to the present case of Schr\"odinger
operators on the line. In particular, if one imposes the additional exponential
falloff of the potential $V$ of the type $V\in L^1(\bbR; \exp(a|x|)dx)$ for
some $a>0$, then $\cF$ and hence the Fredholm determinant on the left-hand side
of \eqref{4.58} permit an analytic continuation through the essential
spectrum of $H$ into a strip of width $a/2$ (w.r.t. the variable $z^{1/2}$).
This is of relevance to the study of resonances of $H$ (cf., e.g.,
\cite{Fr97}, \cite{Si00}, and the literature cited therein).

The result \eqref{4.58} is well-known (although, typically under the
additional assumption that $V$ be real-valued), see, for instance,
\cite{Ge86}, 
\cite[Appendix A]{Ne80}, \cite[Proposition 5.7]{Si79}, \cite{Si00}. Again, the
derivation just presented appears to be the most streamlined available for the
reasons outlined after Remark \ref{r4.3a}.

For an explicit expansion of Fredholm determinants of the type \eqref{4.16} and
\eqref{4.48} (valid in the case of general Green's functions $G$ of
Schr\"odinger operators $H$, not just for $G^{(0)}$ associated with $H^{(0)}$)
we refer to Proposition 2.8 in \cite{Si77}.  

\medskip

Next, we revisit the result \eqref{4.58} from a different and perhaps
somewhat unusual perspective. We intend to rederive the analogous result in
the context of $2$-modified determinants $\det_2(\cdot)$ by rewriting
the scalar second-order Schr\"odinger equation as a first-order
$2\times 2$ system, taking the latter as our point of departure.

Assuming hypothesis \ref{4.39} for the rest of this example, the Schr\"odinger
equation
\begin{equation}
-\psi''(z,x) + V(x)\psi(z,x)=z\psi(z,x), \lb{4.100}
\end{equation} 
is equivalent to the first-order system 
\begin{equation}
\Psi'(z,x)=\begin{pmatrix} 0 & 1 \\ V(x)-z & 0 \end{pmatrix} \Psi(z,x), \quad 
\Psi(z,x) = \begin{pmatrix} \psi(z,x) \\ \psi'(z,x) \end{pmatrix}. \lb{4.101}
\end{equation}
Since $\Phi^{(0)}$ defined by 
\begin{equation}
\Phi^{(0)} (z,x)=\begin{pmatrix} \exp(-iz^{1/2}x) & \exp(iz^{1/2}x) \\
-iz^{1/2} \exp(-iz^{1/2}x) & iz^{1/2} \exp(iz^{1/2}x) \end{pmatrix}, \quad 
\Im(z^{1/2})\geq 0 \lb{4.102}
\end{equation}
with
\begin{equation}
{\det}_{\bbC^2}(\Phi^{(0)} (z,x))=1, \quad (z,x)\in\bbC\times\bbR, \lb{4.103}
\end{equation}
is a fundamental matrix of the system \eqref{4.101} in the case $V=0$ a.e., and
since 
\begin{equation}
\Phi^{(0)} (z,x) \Phi^{(0)} (z,x')^{-1}=\begin{pmatrix} 
\cos(z^{1/2}(x-x')) & z^{-1/2}\sin(z^{1/2}(x-x')) \\
-z^{1/2}\sin(z^{1/2}(x-x'))  &  \cos(z^{1/2}(x-x')) \end{pmatrix}, \lb{4.104}
\end{equation} 
the system \eqref{4.101} has the following pair of linearly independent
solutions for $z\neq 0$,
\begin{align}
&F_\pm(z,x)=F^{(0)}_\pm (z,x) \no \\
& \quad - \int_x^{\pm\infty} dx' \begin{pmatrix} 
\cos(z^{1/2}(x-x')) & z^{-1/2}\sin(z^{1/2}(x-x')) \\
-z^{1/2}\sin(z^{1/2}(x-x'))  &  \cos(z^{1/2}(x-x')) \end{pmatrix} \no \\
& \hspace{2.3cm} \times \begin{pmatrix} 0 & 0 \\ V(x') & 0 \end{pmatrix}
F_\pm(z,x') \no \\ 
& \quad =F^{(0)}_\pm (z,x) - \int_x^{\pm\infty} dx' 
\begin{pmatrix} z^{-1/2}\sin(z^{1/2}(x-x')) & 0 \\
\cos(z^{1/2}(x-x')) & 0 \end{pmatrix} V(x') F_\pm(z,x'), \lb{4.105}\\
& \hspace*{6.8cm} \Im(z^{1/2})\geq 0, \; z \neq 0, \; x\in\bbR, \no
\end{align}
where we abbreviated
\begin{equation}
F^{(0)}_\pm (z,x) = \begin{pmatrix} 1 \\ \pm iz^{1/2} \end{pmatrix} 
\exp(\pm i z^{1/2} x). \lb{4.106}
\end{equation}
By inspection, the first component of \eqref{4.105} is equivalent to
\eqref{4.42} and the second component to the $x$-derivative of \eqref{4.42},
that is, one has
\begin{equation}
F_\pm (z,,x)=\begin{pmatrix} f_\pm (z,x) \\ f'_\pm(z,x) \end{pmatrix}, \quad 
\Im(z^{1/2})\geq 0, \; z \neq 0, \; x\in\bbR. \lb{4.107}
\end{equation}
Next, one introduces 
\begin{align}
\begin{split}
f_1(z,x)&=-u(x) \begin{pmatrix} 1 \\ i z^{1/2} \end{pmatrix} 
\exp( i z^{1/2} x),  \\ 
f_2(z,x)&= -u(x) \begin{pmatrix} 1 \\ -i z^{1/2}
\end{pmatrix} \exp( -i z^{1/2} x), \lb{4.108} \\
g_1(z,x)&=v(x)\bigg(\f{i}{2z^{1/2}} \exp(-iz^{1/2} x) \quad 0\bigg), \\
g_2(z,x)&=v(x)\bigg(\f{i}{2z^{1/2}} \exp(iz^{1/2} x) \quad 0\bigg)
\end{split}
\end{align}
and hence
\begin{align}
H(z,x,x')& = f_1(z,x) g_1(z,x')-f_2(z,x)g_2(z,x') \lb{4.109} \\
& = u(x) \begin{pmatrix} z^{-1/2}\sin(z^{1/2}(x-x')) & 0 \\
\cos(z^{1/2}(x-x')) & 0 \end{pmatrix} v(x') \lb{4.110}
\end{align}
and we introduce
\begin{align}
\widetilde K(z,x,x')&=\begin{cases} f_1(z,x)g_1(z,x'), & x'<x, \\ 
f_2(z,x)g_2(z,x'), & x<x',\end{cases}  \lb{4.111} \\
&=\begin{cases} -u(x)\f{1}{2}\exp(iz^{1/2}(x-x'))\begin{pmatrix} iz^{-1/2} & 0
\\  -1 & 0  \end{pmatrix} v(x'), & x'<x, \\ 
-u(x)\f{1}{2}\exp(-iz^{1/2}(x-x'))\begin{pmatrix} iz^{-1/2} & 0 \\  
1 & 0  \end{pmatrix} v(x'), & x<x', \end{cases} \lb{4.112} \\
& \hspace*{4.8cm} \Im(z^{1/2})\geq 0, \, z\neq 0, \; x,x'\in\bbR. \no  
\end{align}
We note that $\widetilde K(z,\cdot,\cdot)$ is discontinuous on the diagonal
$x=x'$.  Since 
\begin{equation}
\widetilde K(z,\cdot,\cdot)\in L^2(\bbR^2;dx\,dx'), \quad \Im(z^{1/2})\geq 0,
\, z\neq 0,  \lb{4.113}
\end{equation}
the associated operator $\widetilde K(z)$ with integral kernel \eqref{4.112} is
Hilbert--Schmidt,
\begin{equation}
\widetilde K(z) \in \cB_2 (L^2(\bbR;dx)), \quad \Im(z^{1/2})\geq 0, \; 
z\neq 0.  \lb{4.114}
\end{equation}
Next, assuming temporarily that 
\begin{equation}
\supp(V) \, \text{ is compact,} \lb{4.115}
\end{equation}
the integral equations defining $\hat f_j(z,x)$, $j=1,2$,  
\begin{align}
\hat f_1(z,x)&=f_1(z,x)-\int_x^\infty dx'\, H(z,x,x')\hat f_1(z,x'), 
\lb{4.116} \\ 
\hat f_2(z,x)&=f_2(z,x)+\int_{-\infty}^x dx'\, H(z,x,x')\hat f_2(z,x'),
\lb{4.117} \\ 
&\hspace*{2cm}  \Im(z^{1/2})\geq 0, \; z\neq 0, \; x \in\bbR, \no
\end{align}  
yield solutions $\hat f_j(z,\cdot) \in L^2(\bbR; dx)$, $j=1,2$. By
comparison with \eqref{4.105}, one then identifies
\begin{align}
\hat f_1(z,x)&=-u(x) F_+ (z,x), \lb{4.118} \\
\hat f_2(z,x)&=-u(x) F_- (z,x). \lb{4.119}
\end{align} 
We note that the temporary compact support assumption \eqref{4.115} on $V$ has
only been introduced to guarantee that $f_j(z,\cdot), \hat f_j(z,\cdot) \in
L^2(\bbR; dx)^2$, $j=1,2$. This extra hypothesis will soon be removed. 

An application of Lemma \ref{l2.6} and Theorem \ref{t3.3} then yields the
following result.

\begin{theorem} \lb{t4.5a}
Suppose $V\in L^1(\bbR; dx)$ and let 
$z\in\bbC$ with $\Im(z^{1/2})\geq 0$, $z\neq 0$. Then 
\begin{align}
{\det}_2(I-\widetilde K(z))&=\cF(z)\exp\bigg(-\f{i}{2z^{1/2}} \int_\bbR dx\,
V(x)\bigg) \lb{4.120} \\ 
&={\det}_2(I-K(z))  \lb{4.121}
\end{align} 
with $K(z)$ defined in \eqref{4.47}. 
\end{theorem}
\begin{proof}
Assuming temporarily that $\supp(V)$ is compact (cf.\ \eqref{4.115})
equation \eqref{4.120} directly follows from combining \eqref{3.26} (or
\eqref{3.29}) with $a=-\infty$, $b=\infty$, \eqref{3.12} (or \eqref{3.14}),
\eqref{4.58}, and \eqref{4.108}. Equation \eqref{4.121} then follows from
\eqref{3.23}, \eqref{3.2} (or \eqref{3.3}), and \eqref{4.108}. To extend the
result to general $V\in L^1(\bbR; dx)$ one follows the approximation argument
presented in Theorem \ref{t4.3}.
\end{proof}

One concludes that the scalar second-order equation \eqref{4.100} and the
first-order system \eqref{4.101} share the identical $2$-modified
Fredholm determinant. 

\begin{remark} \lb{r4.5b}
Let $\Im(z^{1/2})\geq 0$, $z\neq 0$, and $x\in\bbR$. Then following up on
Remark \ref{r2.8}, one computes
\begin{align}
A(z,x)&=\begin{pmatrix} g_1(z,x)f_1(z,x) & g_1(z,x)f_2(z,x) \\ 
-g_2(z,x)f_1(z,x) & -g_2(z,x)f_2(z,x) \end{pmatrix} \no \\
&= -\f{i}{2z^{1/2}}V(x) \begin{pmatrix} 1 & e^{-2iz^{1/2}x}
\\ -e^{2iz^{1/2}x} & -1 \end{pmatrix} \lb{4.122} \\
&=-\f{i}{2z^{1/2}}V(x) \begin{pmatrix} e^{-iz^{1/2}x} & 0
\\ 0 & e^{iz^{1/2}x} \end{pmatrix} \begin{pmatrix} 1 & 1
\\ -1 & -1 \end{pmatrix} \begin{pmatrix} e^{iz^{1/2}x} & 0
\\ 0 & e^{-iz^{1/2}x} \end{pmatrix}. \no
\end{align}
Introducing 
\begin{equation}
W(z,x)= e^{M(z)x} U(z,x), \quad M(z)=iz^{1/2}\begin{pmatrix} 1 & 0 \\ 
0 & -1 \end{pmatrix}, \lb{4.123}  
\end{equation}
and recalling 
\begin{equation}
U'(z,x)=A(z,x)U(z,x),   \lb{4.124}
\end{equation}
$($cf.\ \eqref{2.18}$)$, equation \eqref{4.124} reduces to 
\begin{equation}
W'(z,x)=\bigg[iz^{1/2}\begin{pmatrix} 1 & 0 \\ 0 & -1\end{pmatrix}
-\f{i}{2z^{1/2}} V(x) \begin{pmatrix} 1 & 1 \\ -1 & -1
\end{pmatrix} \bigg]W(z,x). \lb{4.125}  
\end{equation}
Moreover, introducing 
\begin{equation}
T(z)=\begin{pmatrix} 1 & 1 \\ iz^{1/2} & -iz^{1/2} \end{pmatrix}, 
\quad \Im(z^{1/2})\geq 0, \, z\neq 0, \lb{4.126} 
\end{equation}
one obtains 
\begin{align}
&\bigg[iz^{1/2}\begin{pmatrix} 1 & 0 \\ 0 & -1\end{pmatrix}
-\f{i}{2z^{1/2}} V(x) \begin{pmatrix} 1 & 1 \\ -1 & -1
\end{pmatrix}\bigg]=T(z)^{-1}\begin{pmatrix} 0 & 1 \\ V(x)-z & 0
\end{pmatrix} T(z), \no \\
& \hspace*{6.5cm} \Im(z^{1/2})\geq 0, \, z\neq 0, \; x\in\bbR,
\lb{4.127} 
\end{align}
which demonstrates the connection between \eqref{2.18}, \eqref{4.125},
and \eqref{4.101}.
\end{remark}

\medskip

Finally, we turn to the case of periodic Schr\"odinger operators of
period $\omega>0$: \\
{\bf The case $\mathbf{(a,b)=(0,\omega)}$:}  Assuming 
\begin{equation}
V\in L^1((0,\omega); dx),  \lb{4.64a}
\end{equation}
we introduce two one-parameter families of closed operators in
$L^2((0,\omega); dx)$ defined by 
\begin{align}
&H^{(0)}_\theta f=-f'', \no \\
&f\in \dom\big(H^{(0)}_\theta\big)=\{g\in L^2((0,\omega); dx) \,|\, g,g'\in
AC ([0,\omega]); \, g(\omega)=e^{i\theta}g(0), \no \\
& \hspace*{2.93cm}  g'(\omega)=e^{i\theta}g'(0), \, 
g''\in L^2((0,\omega); dx)\}, \lb{4.64b} \\ 
&H_\theta f=-f''+Vf, \no \\
&f\in \dom(H_\theta)=\{g\in L^2((0,\omega); dx) \,|\, g,g'\in
AC ([0,\omega]); \, g(\omega)=e^{i\theta}g(0), \no \\
& \hspace*{2.7cm} g'(\omega)=e^{i\theta}g'(0), 
\, (-g''+Vg)\in L^2((0,\omega); dx)\}, \lb{4.65}
\end{align}
where $\theta\in [0,2\pi)$. As in the previous cases considered,
$H^{(0)}_\theta$ is self-adjoint and $H_\theta$ is self-adjoint if and
only if $V$ is real-valued.

Introducing the fundamental system of solutions
$c(z,\cdot)$ and $s(z,\cdot)$ of $-\psi''(z)+V\psi(z)=z\psi(z)$, $z\in\bbC$, by
\begin{equation}
c(z,0)=1=s'(z,0), \quad c'(z,0)=0=s(z,0), \lb{4.66}
\end{equation}
the associated fundamental matrix of solutions $\Phi(z,x)$ is defined by
\begin{equation}
\Phi(z,x)=\begin{pmatrix} c(z,x) & s(z,x) \\ c'(z,x) & s'(z,x) \end{pmatrix}.  
\lb{4.67}
\end{equation}
The monodromy matrix is then given by $\Phi(z,\omega)$, and the Floquet
discriminant $\Delta(z)$ is defined as half of the trace of the latter,
\begin{equation}
\Delta(z)=\tr_{\bbC^2}(\Phi(z,\omega))/2=[c(z,\omega)+s'(z,\omega)]/2.
\lb{4.68}
\end{equation}
Thus, the eigenvalue equation for $H_\theta$ reads,
\begin{equation}
\Delta(z)=\cos(\theta). \lb{4.68a}
\end{equation}
In the special case $V=0$ a.e.\ one obtains
\begin{equation}
c^{(0)}(z,x)=\cos(z^{1/2}x), \quad s^{(0)}(z,x)=\sin(z^{1/2}x) \lb{4.68b}
\end{equation}
and hence,
\begin{equation}
\Delta^{(0)}(z)=\cos(z^{1/2}\omega). \lb{4.68c}
\end{equation}

Next we introduce additional solutions $\varphi_\pm (z,\cdot)$, $\psi_\pm
(z,\cdot)$ of $-\psi''(z)+V\psi(z)=z\psi(z)$, $z\in\bbC$, by 
\begin{align}
\varphi_\pm (z,x)&=e^{\pm iz^{1/2}x}+\int_0^{x} dx'
g^{(0)}(z,x,x')V(x')\varphi_\pm (z,x'), \lb{4.69} \\
\psi_\pm (z,x)&=e^{\pm iz^{1/2}x}-\int_x^{\omega} dx'
g^{(0)}(z,x,x')V(x')\psi_\pm (z,x'), \lb{4.70} \\
&\hspace*{3.64cm} \Im(z^{1/2})\geq 0, \; x\in [0,\omega], \no 
\end{align}
where $g^{(0)}(z,x,x')$ is still given by \eqref{4.7}. We also introduce the
Green's function of $H^{(0)}_\theta$,
\begin{align}
G^{(0)}_\theta (z,x,x')&=\big(H^{(0)}_\theta -z\big)^{-1}(x,x') \no \\
&=\f{i}{2z^{1/2}}\bigg[e^{iz^{1/2}|x-x'|}+
\f{e^{iz^{1/2}(x-x')}}{e^{i\theta}e^{-iz^{1/2}\omega}-1}
+ \f{e^{-iz^{1/2}(x-x')}}{e^{-i\theta}e^{-iz^{1/2}\omega}-1} \bigg], 
\lb{4.71} \\
& \hspace*{4.9cm} \Im(z^{1/2})>0, \; x,x'\in (0,\omega). \no 
\end{align}
Introducing again the factorization \eqref{4.13} of $V=uv$, one verifies as in
\eqref{4.14} that 
\begin{align}
(H_\theta-z)^{-1}&=\big(H^{(0)}_\theta-z\big)^{-1} \no \\
&\quad -\big(H^{(0)}_\theta-z\big)^{-1}v
\Big[I+\ol{u\big(H^{(0)}_\theta -z\big)^{-1}v}\Big]^{-1} 
u\big(H^{(0)}_\theta -z\big)^{-1}, \lb{4.72} \\ 
& \hspace*{4.15cm}
z\in\bbC\backslash\{\spec(H_\theta)\cup\spec(H_\theta^{(0)})\}. \no  
\end{align}
To establish the connection with the notation used in Sections \ref{s2}
and \ref{s3}, we introduce the operator $K_\theta (z)$ in
$L^2((0,\omega); dx)$  (cf.\ \eqref{2.3}, \eqref{4.15}) by
\begin{equation}
K_{\theta}(z)=-\ol{u\big(H^{(0)}_\theta -z\big)^{-1}v}, \quad
z\in\bbC\backslash\spec\big(H^{(0)}_\theta\big)
\lb{4.73}
\end{equation}
with integral kernel
\begin{equation}
K_\theta(z,x,x')=-u(x)G^{(0)}_{\theta}(z,x,x')v(x'), \quad 
z\in\bbC\backslash\spec\big(H^{(0)}_\theta\big), \; x,x' \in [0,\omega],
\lb{4.74} 
\end{equation}
and the Volterra operators $H_{0} (z)$, $H_\omega (z)$ (cf.\ \eqref{2.4},
\eqref{2.5}) with integral kernel
\begin{equation}
H(z,x,x')=u(x)g^{(0)}(z,x,x')v(x').   \lb{4.75} 
\end{equation}
Moreover, we introduce for a.e.\ $x\in (0,\omega)$, 
\begin{align}
f_1(z,x)&=f_2(z,x)=f(z,x)=-u(x) (e^{iz^{1/2}x} \;\, e^{-iz^{1/2}x}), \no \\
g_1(z,x)&=\f{i}{2z^{1/2}}v(x)\begin{pmatrix} 
\f{\exp(i\theta)\exp(-iz^{1/2}\omega)\exp(-iz^{1/2}x)}{\exp(i\theta)
\exp(-iz^{1/2}\omega)-1} \\
\f{\exp(iz^{1/2}x)}{\exp(-i\theta)\exp(-iz^{1/2}\omega)-1} \end{pmatrix},
\lb{4.76} \\ 
g_2(z,x)&=\f{i}{2z^{1/2}}v(x)\begin{pmatrix} 
\f{\exp(-iz^{1/2}x)}{\exp(i\theta)\exp(-iz^{1/2}\omega)-1} \\
\f{\exp(-i\theta)\exp(-iz^{1/2}\omega)\exp(iz^{1/2}x)}{\exp(-i\theta)
\exp(-iz^{1/2}\omega)-1} \end{pmatrix}. \no 
\end{align}
Introducing $\hat f_j(z,x)$, $j=1,2$, by
\begin{align}
\hat f_1(z,x)&=f(z,x)-\int_x^\omega dx'\, H(z,x,x')\hat f_1(z,x'), \lb{4.77}
\\ 
\hat f_2(z,x)&=f(z,x)+\int_{0}^x dx'\, H(z,x,x')\hat f_2(z,x'),
\lb{4.78} \\ 
&\hspace*{1.75cm}  \Im(z^{1/2})\geq 0, \; z\neq 0, \; x \geq 0, \no
\end{align}  
yields solutions $\hat f_j(z,\cdot) \in L^2((0,\omega); dx)$, $j=1,2$. By
comparison with \eqref{4.5}, \eqref{4.6}, one then identifies
\begin{align}
\hat f_1(z,x)&=-u(x) (\psi_+ (z,x) \;\, \psi_-(z,x)), \lb{4.79} \\
\hat f_2(z,x)&=-u(x) (\varphi_+ (z,x) \;\, \varphi_-(z,x)). \lb{4.80} 
\end{align}

Next we mention the following result.

\begin{theorem} \lb{t4.6}
Suppose $V\in L^1((0,\omega); dx)$, let $\theta\in [0,2\pi)$, and
$z\in\bbC\backslash\spec\big(H^{(0)}_\theta\big)$. Then 
\begin{equation}
K_\theta(z)\in \cB_1(L^2((0,\omega); dx)) \lb{4.81}
\end{equation}
and 
\begin{equation}
\det(I-K_\theta (z))=\f{\Delta(z)-\cos(\theta)}{\cos(z^{1/2}\omega) 
-\cos(\theta)}. \lb{4.82}
\end{equation} 
\end{theorem}
\begin{proof} Since the integral kernel of $K_\theta(z)$ is square integrable
over $(0,\omega)\times (0,\omega)$ one has of course 
$K_\theta(z)\in \cB_2(L^2((0,\omega);dx))$. To prove its trace class property
one imbeds $(0,\omega)$ into $\bbR$ in analogy to the half-line case
discussed in the proof of Theorem \ref{t4.1}, introducing 
\begin{equation}
L^2(\bbR;dx)=L^2((0,\omega); dx)\oplus L^2(\bbR\backslash [0,\omega];dx) 
\lb{4.82a}
\end{equation}
and 
\begin{align}
\begin{split}
\tilde u(x)&=\begin{cases} u(x), & x\in (0,\omega), \\ 0, & x\notin (0,\omega),
\end{cases} \quad 
\tilde v(x)=\begin{cases} v(x), & x\in (0,\omega), \\ 0, & x\notin (0,\omega),
\end{cases} \\\
\widetilde V(x)&=\begin{cases} V(x), & x\in (0,\omega), \\ 0, & x\notin
(0,\omega). \end{cases} \lb{4.82b}
\end{split}
\end{align} 
At this point one can follow the proof of Theorem \ref{t4.1} line by line using
\eqref{4.71} instead of \eqref{4.24cb} and noticing that the second and third
term on the right-hand side of \eqref{4.71} generate rank one terms upon
multiplying them by $\tilde u(x)$ from the left and $\tilde v(x')$ from the
right. 

By \eqref{4.68a} and \eqref{4.68c}, and since 
\begin{equation}
\det(I-K_\theta (z))=\det\Big(\big(H^{(0)}_\theta -z\big)^{-1/2}(H_\theta -z)
\big(H^{(0)}_\theta -z\big)^{-1/2}\Big), \lb{4.83}
\end{equation}
$\det(I-K_\theta (z))$ and $[\Delta(z)-\cos(\theta)]/[\cos(z^{1/2}\omega) 
-\cos(\theta)]$ have the same set of zeros and poles. Moreover, since
either expression satisfies the asymptotics $1+\oh(1)$ as
$z\downarrow -\infty$, one obtains \eqref{4.82}.
\end{proof}

An application of Lemma \ref{l2.6} and Theorem \ref{t3.2} then yields the
following result relating the Fredholm determinant of $I-K_\theta (z)$ and the
Floquet discriminant $\Delta(z)$.

\begin{theorem} \lb{t4.7}
Suppose $V\in L^1((0,\omega); dx)$, let $\theta\in [0,2\pi)$, and
$z\in\bbC\backslash\spec\big(H^{(0)}_\theta\big)$. Then 
\begin{align}
&\det(I-K_\theta (z))=\f{\Delta(z)-\cos(\theta)}{\cos(z^{1/2}\omega) 
-\cos(\theta)} \no \\
&\quad =\bigg[1+\f{i}{2z^{1/2}}
\f{e^{i\theta}e^{-iz^{1/2}\omega}}{e^{i\theta}e^{-iz^{1/2}\omega}-1}
\int_0^\omega dx\, e^{-iz^{1/2}x}V(x)\psi_+(z,x)\bigg] \no \\
&\qquad \times \bigg[1+\f{i}{2z^{1/2}}
\f{1}{e^{-i\theta}e^{-iz^{1/2}\omega}-1}
\int_0^\omega dx\, e^{iz^{1/2}x}V(x)\psi_-(z,x)\bigg] \no \\
&\qquad
+\f{1}{4z}\f{e^{i\theta}e^{-iz^{1/2}\omega}}{\big[
e^{i\theta}e^{-iz^{1/2}\omega}-1\big]
\big[e^{-i\theta}e^{-iz^{1/2}\omega}-1\big]}
\int_0^\omega dx\, e^{iz^{1/2}x}V(x)\psi_+(z,x) \no \\
& \hspace*{6cm} 
\times \int_0^\omega dx\, e^{-iz^{1/2}x}V(x)\psi_-(z,x) \lb{4.84} 
\end{align}
\begin{align}
&\quad \;\, =\bigg[1+\f{i}{2z^{1/2}}
\f{1}{e^{i\theta}e^{-iz^{1/2}\omega}-1}
\int_0^\omega dx\, e^{-iz^{1/2}x}V(x)\varphi_+(z,x)\bigg] \no \\
&\qquad \;\, \times \bigg[1+\f{i}{2z^{1/2}}
\f{e^{-i\theta}e^{-iz^{1/2}\omega}}{e^{-i\theta}e^{-iz^{1/2}\omega}-1}
\int_0^\omega dx\, e^{iz^{1/2}x}V(x)\varphi_-(z,x)\bigg] \no \\
&\qquad \;\,
+\f{1}{4z}\f{e^{-i\theta}e^{-iz^{1/2}\omega}}{\big[
e^{i\theta}e^{-iz^{1/2}\omega}-1\big]
\big[e^{-i\theta}e^{-iz^{1/2}\omega}-1\big]}
\int_0^\omega dx\, e^{iz^{1/2}x}V(x)\varphi_+(z,x) \no \\
& \hspace*{6cm} \;\,
\times \int_0^\omega dx\, e^{-iz^{1/2}x}V(x)\varphi_-(z,x). \lb{4.85}
\end{align} 
\end{theorem}
\begin{proof}
Again Lemma \ref{l2.6}  applies and one infers from \eqref{2.37} and
\eqref{4.76}--\eqref{4.80} that
\begin{equation}
U(z,x)=\begin{pmatrix} 1- \int_x^\omega dx'\, g_1(z,x')\hat
f(z,x') & \int_0^x dx'\, g_1(z,x')\hat f(z,x') \\ 
\int_x^\omega dx'\, g_2(z,x')\hat f(z,x') & 1-\int_0^x
dx'\, g_2(z,x')\hat f(z,x') \end{pmatrix}, \quad x\in [0,\omega], \lb{4.86}
\end{equation}
becomes
\begin{align}
U_{1,1}(z,x)&=I_2+ \f{i}{2z^{1/2}}\int_x^\omega dx' 
\begin{pmatrix} \f{\exp(i\theta)\exp(-iz^{1/2}\omega)
\exp(-iz^{1/2}x')}{\exp(i\theta)\exp(-iz^{1/2}\omega)-1} \\
\f{e^{iz^{1/2}x'}}{\exp(-i\theta)\exp(-iz^{1/2}\omega)-1}\end{pmatrix}
V(x') \no \\
& \hspace*{5.75cm} \times (\psi_+ (z,x') \;\, \psi_-(z,x')), \lb{4.87} \\
U_{1,2}(z,x)&= -\f{i}{2z^{1/2}}\int_0^x dx' 
\begin{pmatrix} \f{\exp(i\theta)\exp(-iz^{1/2}\omega)
\exp(-iz^{1/2}x')}{\exp(i\theta)\exp(-iz^{1/2}\omega)-1} \\
\f{\exp(iz^{1/2}x')}{\exp(-i\theta)\exp(-iz^{1/2}\omega)-1}\end{pmatrix}
V(x') \no \\
& \hspace*{5.2cm} \times (\varphi_+ (z,x') \;\, \varphi_-(z,x')), \lb{4.88} \\
U_{2,1}(z,x)&=-\f{i}{2z^{1/2}}\int_x^\omega dx' 
\begin{pmatrix} \f{\exp(-iz^{1/2}x')}{\exp(i\theta)\exp(-iz^{1/2}\omega)-1} \\
\f{\exp(-i\theta)\exp(-iz^{1/2}\omega)\exp(iz^{1/2}x')}{\exp(-i\theta)
\exp(-iz^{1/2}\omega)-1}\end{pmatrix}
V(x') \no \\
& \hspace*{5.25cm} \times (\psi_+ (z,x') \;\, \psi_-(z,x')), \lb{4.89} \\
U_{2,2}(z,x)&=I_2+\f{i}{2z^{1/2}}\int_0^x dx' 
\begin{pmatrix} \f{\exp(-iz^{1/2}x')}{\exp(i\theta)\exp(-iz^{1/2}\omega)-1} \\
\f{\exp(-i\theta)\exp(-iz^{1/2}\omega)\exp(iz^{1/2}x')}{\exp(-i\theta)
\exp(-iz^{1/2}\omega)-1}\end{pmatrix}
V(x') \no \\
& \hspace*{5.75cm} \times (\varphi_+ (z,x') \;\, \varphi_-(z,x')). \lb{4.90}
\end{align}
Relations \eqref{3.5} and \eqref{3.7} of Theorem \ref{t3.2} with $m=1$,
$n_1=n_2=2$, $n=4$, then immediately yield \eqref{4.84} and \eqref{4.85}.  
\end{proof} 

To the best of our knowledge, the representations \eqref{4.84} and
\eqref{4.85} of $\Delta(z)$ appear to be new. They are the analogs of the
well-known representations of Jost functions \eqref{4.10}, \eqref{4.11} and
\eqref{4.45} on the half-line and on the real line, respectively. That 
the Floquet discriminant $\Delta(z)$ is related to infinite determinants is
well-known. However, the connection between $\Delta(z)$ and determinants of
Hill-type discussed in the literature (cf., e.g., \cite{Ma55},   
\cite[Ch.\ III, Sect.\ VI.2]{GGK00}, \cite[Sect.\ 2.3]{MW79}) is of a different
nature than the one in \eqref{4.82} and based on the Fourier expansion of the
potential $V$. For different connections between Floquet theory and
perturbation determinants we refer to \cite{GW95}.

\section{Integral operators of convolution-type with rational symbols}
\lb{s5}

In our final section we rederive the explicit formula for the
$2$-modified Fredholm determinant corresponding to integral operators
of  convolution-type, whose integral kernel is associated with a symbol
given by a rational function, in an elementary and straghtforward
manner. This determinant formula represents a truncated Wiener--Hopf
analog of Day's formula for the determinant associated with finite
Toeplitz matrices generated by the Laurent expansion of a rational
function.

Let $\tau>0$. We are interested in truncated Wiener--Hopf-type operators
$K$ in $L^2((0,\tau); dx)$ of the form
\begin{equation}
(Kf)(x) = \int_0^\tau dx'\, k(x-x') f(x'), \quad f\in L^2((0,\tau);dx),
\lb{5.1}
\end{equation}
where $k(\cdot)$, extended from $[-\tau,\tau]$ to $\bbR\backslash\{0\}$,
is defined by
\begin{equation}
k(t)=\begin{cases} \sum_{\ell\in \cL} \alpha_\ell \, e^{-\lambda_\ell
t},  & t>0, \\
\sum_{m\in \cM} \beta_m e^{\mu_m t}, & t<0 \end{cases} \lb{5.2}
\end{equation}
and 
\begin{align}
\begin{split}
&\alpha_\ell\in\bbC, \; \ell\in\cL=\{1,\dots,L\}, \; L\in\bbN, \\
&\beta_m\in\bbC, \; m\in\cM=\{1,\dots,M\}, \; M\in\bbN, \\
&\lambda_\ell\in\bbC, \; \Re(\lambda_\ell)>0, \; \ell\in\cL, \\
&\mu_m\in\bbC, \; \Re(\mu_m)>0, \; m\in\cM. \lb{5.3}
\end{split}
\end{align}
In terms of semi-separable integral kernels, $k$ can be rewritten as,
\begin{equation}
k(x-x')=K(x,x')=\begin{cases} f_1(x)g_1(x'), & 0<x'< x< \tau, \\ 
f_2(x)g_2(x'), & 0<x<x'<\tau, \end{cases}  \lb{5.4}
\end{equation}
where
\begin{align}
\begin{split}
f_1(x)&=\big(\alpha_1  e^{-\lambda_1 x},\dots,\alpha_L 
e^{-\lambda_L x}\big), \\ 
f_2(x)&=\big(\beta_1 e^{\mu_1 x},\dots,\beta_M e^{\mu_M x}\big), \\ 
g_1(x)&=\big(e^{\lambda_1 x},\dots,e^{\lambda_L x}\big)^\top, \\ 
g_2(x)&=\big(e^{-\mu_1 x},\dots,e^{-\mu_M x}\big)^\top. \lb{5.5}
\end{split}
\end{align}
Since $K(\cdot,\cdot)\in L^2((0,\tau)\times (0,\tau); dx\,dx')$, the
operator $K$ in \eqref{5.1} belongs to the Hilbert--Schmidt
class,
\begin{equation}
K\in\cB_2(L^2((0,\tau); dx)). \lb{5.5a}
\end{equation}  

Associated with $K$ we also introduce the Volterra operators
$H_0$, $H_\tau$ (cf.\ \eqref{2.4}, \eqref{2.5}) in $L^2((0,\tau);dx)$
with integral kernel
\begin{equation}
h(x-x')=H(x,x')=f_1(x)g_1(x')-f_2(x)g_2(x'), \lb{5.6}
\end{equation}
such that
\begin{equation}
h(t)=\sum_{\ell\in \cL} \alpha_\ell \, e^{-\lambda_\ell t}-
\sum_{m\in \cM} \beta_m e^{\mu_m t}. \lb{5.7}
\end{equation}
In addition, we introduce the Volterra integral equation
\begin{equation}
\hat f_2(x)=f_2(x)+\int_{0}^x dx'\, h(x-x')\hat f_2(x'), \quad 
 x \in(0,\tau) \lb{5.8}
\end{equation}
with solution $\hat f_2\in L^2((0,\tau);dx)$. 

Next, we introduce the Laplace transform $\bbF$ of a function $f$ by
\begin{equation}
\bbF(\zeta)=\int_0^\infty dt\, e^{-\zeta t} f(t), \lb{5.9}
\end{equation}
where either $f\in L^r((0,\infty);dt)$, $r\in\{1,2\}$ and $\Re(\zeta)>0$,
or, $f$ satisfies an exponential bound of the type $|f(t)|\leq C\exp(Dt)$
for some $C>0$, $D\geq 0$ and then $\Re(\zeta)>D$. Moreover, whenever
possible, we subsequently meromorphically continue $\bbF$ into the
half-plane $\Re(\zeta)<0$ and $\Re(\zeta)<D$, respectively, and for
simplicity denote the result again by $\bbF$.

Taking the Laplace transform of equation \eqref{5.8}, one obtains
\begin{equation}
\hatt \bbF_2 (\zeta)=\bbF_2(\zeta)+\bbH(\zeta)\hatt \bbF_2 (\zeta),  
\lb{5.10}
\end{equation}
where
\begin{align}
\bbF_2(\zeta)&=\big(\beta_1(\zeta-\mu_1)^{-1},\dots,
\beta_M(\zeta-\mu_M)^{-1}\big), \lb{5.11} \\
\bbH(\zeta)&=\sum_{\ell\in\cL} \alpha_\ell (\zeta+\lambda_\ell)^{-1} 
-\sum_{m\in\cM} \beta_m (\zeta-\mu_m)^{-1} \lb{5.12}
\end{align}
and hence solving \eqref{5.10}, yields
\begin{equation}
\hatt \bbF_2 (\zeta)=(1-\bbH(\zeta))^{-1}\big(\beta_1
(\zeta-\mu_1)^{-1},\dots,\beta_M(\zeta-\mu_M)^{-1}\big). \lb{5.13}
\end{equation}
Introducing the Fourier transform $\cF(k)$ of the kernel function $k$ by
\begin{equation}
\cF(k)(x)=\int_\bbR dt\, e^{ixt} k(t), \quad x\in\bbR, \lb{5.14}
\end{equation}
one obtains the rational symbol
\begin{equation}
\cF(k)(x)=\sum_{\ell\in\cL} \alpha_\ell (\lambda_\ell-ix)^{-1} + 
\sum_{m\in\cM} \beta_m (\mu_m+ix)^{-1}. \lb{5.15}
\end{equation}
Thus, 
\begin{equation}
1-\bbH(-ix)=1-\cF(k)(x)= \prod_{n\in\cN} (-ix+i\zeta_n)
\prod_{\ell\in\cL} (-ix+\lambda_\ell)^{-1} 
\prod_{m\in\cM} (-ix-\mu_m)^{-1} \lb{5.16}
\end{equation}
for some
\begin{equation}
\zeta_n\in\bbC, \; n\in\cN=\{1,\dots,N\}, \; N=L+M. \lb{5.17}
\end{equation}
Consequently,
\begin{align}
1-\bbH(\zeta)&= \prod_{n\in\cN} (\zeta+i\zeta_n)
\prod_{\ell\in\cL} (\zeta+\lambda_\ell)^{-1} 
\prod_{m\in\cM} (\zeta-\mu_m)^{-1}, 
\lb{5.18} \\
(1-\bbH(\zeta))^{-1}&=1+\sum_{n\in\cN} \gamma_n (\zeta+i\zeta_n)^{-1}, 
\lb{5.19}
\end{align}
where
\begin{equation}
\gamma_n=\prod_{\substack{n'\in\cN\\n'\neq n}}
(i\zeta_n-i\zeta_{n'})^{-1} \prod_{\ell\in\cL} (\lambda_\ell-i\zeta_n)
\prod_{m\in\cM} (-i\zeta_n-\mu_m), \quad n\in\cN. \lb{5.20}
\end{equation}
Moreover, one computes
\begin{equation}
\beta_m=\prod_{\ell\in\cL} (\mu_m+\lambda_\ell)^{-1} 
\prod_{\substack{m'\in\cM\\m'\neq m}} (\mu_m-\mu_{m'})^{-1} 
\prod_{n\in\cN} (\mu_m+i\zeta_n), \quad m\in\cM. \lb{5.21}
\end{equation}
Combining \eqref{5.13} and \eqref{5.19} yields 
\begin{equation}
\hatt\bbF_2(\zeta)=\bigg(1+\sum_{n\in\cN}\gamma_n(\zeta+i\zeta_n)^{-1}
\bigg)\big(\beta_1(\zeta-\mu_1)^{-1},\dots,\beta_M(\zeta-\mu_M)^{-1}
\big) \lb{5.22}
\end{equation}
and hence
\begin{align}
\hat f_2(x)&=\bigg(\beta_1\bigg[e^{\mu_1 x}-\sum_{n\in\cN} \gamma_n
\big(e^{-i\zeta_n x}-e^{\mu_1 x}\big)(\mu_1+i\zeta_n)^{-1}\bigg],
\dots \no \\
& \qquad \, \dots,\beta_M\bigg[e^{\mu_M x}-\sum_{n\in\cN} \gamma_n 
\big(e^{-i\zeta_n x}-e^{\mu_M x}\big)(\mu_M+i\zeta_n)^{-1}\bigg] 
\bigg). \lb{5.23}
\end{align}
In view of \eqref{3.29} we now introduce the $M\times M$ matrix
\begin{equation}
G=\big(G_{m,m'}\big)_{m,m'\in\cM}=\int_0^\tau dx\, g_2(x)\hat f_2(x).
\lb{5.24}
\end{equation}
\begin{lemma} \lb{l5.1}
One computes
\begin{equation}
G_{m,m'}=\delta_{m,m'}+e^{-\mu_m \tau}\beta_{m'}\sum_{n\in\cN}\gamma_n 
e^{-i\zeta_n \tau}(\mu_m+i\zeta_n)^{-1}(\mu_{m'}+i\zeta_n)^{-1}, 
\quad m,m'\in\cM. \lb{5.25}
\end{equation}
\end{lemma}
\begin{proof}
By \eqref{5.24},
\begin{align}
G_{m,m'}&=\int_0^\tau dt\, e^{-\mu_m t}\beta_{m'}\bigg(e^{\mu_{m'}t}
-\sum_{n\in\cN} \gamma_n \big(e^{-i\zeta_n t}-e^{\mu_{m'}t}
\big)(i\zeta_n+\mu_{m'})^{-1}\bigg) \no \\
&=\beta_{m'}\int_0^\tau dt\, e^{-(\mu_m-\mu_{m'})t}\bigg(1+
\sum_{n\in\cN} \gamma_n (i\zeta_n+\mu_{m'})^{-1}\bigg) \no \\
& \quad -\beta_{m'}\int_0^\tau dt\, e^{-\mu_m t}
\sum_{n\in\cN} \gamma_n e^{-i\zeta_n t} (i\zeta_n+\mu_{m'})^{-1} \no \\
&= -\beta_{m'}\sum_{n\in\cN} \gamma_n (i\zeta_n+\mu_{m'})^{-1}
\int_0^\tau dt\, e^{-(i\zeta_n+\mu_m) t} \no \\
&=\beta_{m'}\sum_{n\in\cN} \gamma_n
\big[e^{-(i\zeta_n+\mu_m)t}-1\big](i\zeta_n+\mu_m)^{-1}
(i\zeta_n+\mu_{m'})^{-1}. \lb{5.26}
\end{align}
Here we used the fact that
\begin{equation}
1+\sum_{n\in\cN} \gamma_n (i\zeta_n+\mu_{m'})^{-1}=0, \lb{5.27}
\end{equation}
which follows from
\begin{equation}
1+\sum_{n\in\cN} \gamma_n
(i\zeta_n+\mu_{m'})^{-1}=(1-\bbH(\mu_{m'}))^{-1}=0, \lb{5.28}
\end{equation}
using \eqref{5.18} and \eqref{5.19}. Next, we claim that
\begin{equation}
-\beta_{m'}\sum_{n\in\cN} \gamma_n (i\zeta_n+\mu_m)^{-1}
(i\zeta_n+\mu_{m'})^{-1}=\delta_{m,m'}. \lb{5.29}
\end{equation}
Indeed, if $m\neq m'$, then
\begin{align}
&\sum_{n\in\cN}\gamma_n (i\zeta_n+\mu_m)^{-1}(i\zeta_n+\mu_{m'})^{-1} 
\no \\
&\quad =-\sum_{n\in\cN}\gamma_n (\mu_m-\mu_{m'})^{-1}
\big[(i\zeta_n+\mu_m)^{-1}
-(i\zeta_n+\mu_{m'})^{-1}\big]=0, \lb{5.30}
\end{align}
using \eqref{5.27}. On the other hand, if $m=m'$, then
\begin{align}
\beta_{m}\sum_{n\in\cN} \gamma_n (i\zeta_n+\mu_{m})^{-2}
&=-\beta_{m}\f{d}{d\zeta} (1-\bbH(\zeta))^{-1}\bigg|_{\zeta=\mu_{m}} 
\no \\
&=\underset{\zeta=\mu_{m}}{\text{Res}} \big(\bbH(\zeta)\big)
\f{d}{d\zeta}(1-\bbH(\zeta))^{-1}\bigg|_{\zeta=\mu_{m}} \no \\
&=-\underset{\zeta=\mu_{m}}{\text{Res}} \f{d}{d\zeta} 
\log\big((1-\bbH(\zeta))^{-1}\big) \no \\
&=-1, \lb{5.31}
\end{align}
using \eqref{5.18}. This proves \eqref{5.29}. Combining \eqref{5.26} and
\eqref{5.29} yields \eqref{5.25}.
\end{proof}

Given Lemma \ref{l5.1}, one can decompose $I_M-G$ as
\begin{equation}
I_M-G=\diag(e^{-\mu_1\tau},\dots,e^{-\mu_M\tau})\,\Gamma \,
\diag(\beta_1,\dots,\beta_M), \lb{5.32}
\end{equation}
where $\diag(\cdot)$ denotes a diagonal matrix and the $M\times M$
matrix $\Gamma$ is defined by
\begin{equation}
\Gamma=\big(\Gamma_{m,m'}\big)_{m,m'\in\cM}=\bigg(-\sum_{n\in\cN} 
\gamma_n e^{-i\zeta_n \tau}
(\mu_m+i\zeta_n)^{-1}(\mu_{m'}+i\zeta_n)^{-1}\bigg)_{m,m'\in\cM}. 
\lb{5.33}
\end{equation}
The matrix $\Gamma$ permits the factorization
\begin{equation}
\Gamma=A\,\diag(\gamma_1 e^{-i\zeta_1\tau},\dots,\gamma_N 
e^{-i\zeta_N\tau})\, B, \lb{5.34}
\end{equation}
where $A$ is the $M\times N$ matrix
\begin{equation}
A=\big(A_{m,n}\big)_{m\in\cM,n\in\cN}
=\big((\mu_m+i\zeta_n)^{-1}\big)_{m\in\cM,n\in\cN} \lb{5.35}
\end{equation}
and $B$ is the $N\times M$ matrix
\begin{equation}
B=\big(B_{n,m}\big)_{n\in\cN,m\in\cM}
=\big(-(\mu_{m}+i\zeta_n)^{-1}\big)_{n\in\cN,m\in\cM} \, . \lb{5.36}
\end{equation}

Next, we denote by $\Psi$ the set of all monotone functions $\psi\colon 
\{1,\dots,M\}\to \{1,\dots,N\}$ (we recall $N=L+M$) such that
\begin{equation}
\psi(1)<\dots < \psi(M). \lb{5.37}
\end{equation} 
The set $\Psi$ is in a one-to-one correspondence with all subsets
${\wti \cM}^\bot=\{1,\dots,N\}\backslash\wti \cM$ of $\{1,\dots,N\}$
which consist of $L$ elements. Here $\wti\cM\subseteq\{1,\dots,N\}$ with
cardinality of $\cM$ equal to $M$, $|\wti\cM|=M$. 

Moreover, denoting by $A_\psi$ and $B^\psi$ the $M\times M$ matrices
\begin{align}
A_\psi&=\big(A_{m,\psi(m')}\big)_{m,m'\in\cM}, \quad \psi\in\Psi, 
\lb{5.38} \\ 
B^\psi&=\big(B_{\psi(m),m'}\big)_{m,m'\in\cM}, \quad \psi\in\Psi, 
\lb{5.39}
\end{align}
one notices that
\begin{equation}
A_\psi^\top=-B^\psi, \quad \psi\in\Psi.  \lb{5.40}
\end{equation}
The matrix $A^\psi$ is of Cauchy-type and one infers (cf.\
\cite[p.\ 36]{Kn75}) that
\begin{equation}
A_\psi^{-1}=D_1^\psi A_\psi^\top D_2^\psi, \lb{5.41}
\end{equation}
where $D_j^\psi$, $j=1,2$, are diagonal matrices with diagonal
entries given by
\begin{align}
\big(D_1^\psi\big)_{m,m}&=\prod_{m'\in\cM}(\mu_{m'}+i\zeta_{\psi(m)})
\prod_{\substack{m''\in\cM\\m''\neq m}}
(-i\zeta_{\psi(m'')}+i\zeta_{\psi(m)})^{-1}, \quad m\in\cM, \lb{5.42} \\
\big(D_2^\psi\big)_{m,m}&=\prod_{m'\in\cM}(\mu_{m}+i\zeta_{\psi(m')})
\prod_{\substack{m''\in\cM\\m''\neq m}}
(\mu_m-\mu_{m''})^{-1}, \quad m\in\cM. \lb{5.43}
\end{align}

One then obtains the following result.
\begin{lemma} \lb{l5.2}
The determinant of $I_M-G$ is of the form
\begin{align}
{\det}_{\bbC^M}(I_M-G)&=(-1)^M\exp\bigg(-\tau\sum_{m\in\cM} \mu_m\bigg)
\bigg(\prod_{\ell\in\cL} \beta_\ell\bigg) \sum_{\psi\in\Psi} 
\bigg(\prod_{\ell'\in\cL} \gamma_{\psi(\ell')}\bigg) \no \\
& \quad\times\exp\bigg(-i\tau\sum_{\ell''\in\cL}
\zeta_{\psi(\ell'')}\bigg) \big[{\det}_{\bbC^M}\big(D_1^\psi\big)
{\det}_{\bbC^M}\big(D_2^\psi\big)\big]^{-1}. \lb{5.44}
\end{align}
\end{lemma}
\begin{proof}
Let $\psi\in\Psi$. Then
\begin{align}
{\det}_{\bbC^M}\big(A_\psi\big){\det}_{\bbC^M}\big(B^\psi\big)&=
(-1)^M \big[{\det}_{\bbC^M}\big(A_\psi\big)\big]^2 \no \\
&=(-1)^M \big[{\det}_{\bbC^M}\big(D_1^\psi\big)
{\det}_{\bbC^M}\big(D_2^\psi\big)\big]^{-1}. \lb{5.45}
\end{align}
An application of the Cauchy--Binet formula for determinants yields
\begin{equation}
{\det}_{\bbC^M}(\Gamma)=\sum_{\psi\in\Psi}{\det}_{\bbC^M}
\big(A_\psi\big) {\det}_{\bbC^M} \big(B^\psi\big) \prod_{m\in\cM} 
\gamma_{\psi(m)} e^{-i\tau \zeta_{\psi(m)}}. \lb{5.46}
\end{equation}
Combining \eqref{5.32}, \eqref{5.45}, and \eqref{5.46} then yields 
\eqref{5.44}.
\end{proof}

Applying Theorem \ref{t3.3} then yields the principal result of this
section.

\begin{theorem} \lb{t5.3}
Let $K$ be the Hilbert--Schmidt operator defined in 
\eqref{5.1}--\eqref{5.3}. Then
\begin{align}
{\det}_2 (I-K)&=\exp\bigg(\tau k(0_-)-\tau\sum_{m\in\cM}
\mu_m \bigg) \sum_{\substack{\wti \cL\subseteq 
\{1,\dots,N\}\\ |\wti \cL|=L}} V_{\wti\cL} \, \exp\big(-i\tau
v_{\wti\cL}\big) \lb{5.47} \\
&=\exp\bigg(\tau k(0_+)-\tau\sum_{\ell\in\cL}
\lambda_\ell \bigg) \sum_{\substack{\wti \cM\subseteq 
\{1,\dots,N\}\\ |\wti \cM|=M}} W_{\wti\cM} \, \exp\big(i\tau
w_{\wti\cM}\big). 
\lb{5.48} 
\end{align}
Here $k(0_\pm)=\lim_{\varepsilon\downarrow 0} k(\pm\varepsilon)$, 
$|\cS|$ denotes the cardinality of $\cS\subset\bbN$, and
\begin{align}
V_{\wti\cL}&=\prod_{\ell\in\cL, \, m\in{\wti\cL}^\bot} 
(\lambda_\ell-i\zeta_m)\prod_{\ell'\in\wti\cL, m'\in\cM}
(\mu_{m'}+i\zeta_{\ell'}) \no \\
&\quad \times \prod_{\ell''\in\cL, m''\in\cM}
(\mu_{m''}+\lambda_{\ell''})^{-1} \prod_{\ell'''\in\wti\cL,
m'''\in{\wti\cL}^\bot} (i\zeta_{m'''}-i\zeta_{\ell'''})^{-1}, \lb{5.49}
\\ 
W_{\wti\cM}&=\prod_{\ell\in\cL, \, m\in{\wti\cM}} 
(\lambda_\ell-i\zeta_m)\prod_{\ell'\in{\wti\cM}^\bot, m'\in\cM}
(\mu_{m'}+i\zeta_{\ell'}) \no \\
&\quad \times \prod_{\ell''\in\cL, m''\in\cM}
(\mu_{m''}+\lambda_{\ell''})^{-1} \prod_{\ell'''\in{\wti\cM}^\bot,
m'''\in\wti\cM} (i\zeta_{\ell'''}-i\zeta_{m'''})^{-1}, \lb{5.50} \\
v_{\wti\cL}&=\sum_{m\in{\wti\cL}^\bot} \zeta_{m}, \lb{5.51} \\
w_{\wti\cM}&=\sum_{\ell\in{\wti\cM}^\bot} \zeta_{\ell} \lb{5.52} 
\end{align}
with
\begin{align}
{\wti \cL}^\bot&=\{1,\dots,N\}\backslash\wti\cL \, \text{ for } \,
\wti \cL\subseteq \{1,\dots,N\}, \;  |\wti \cL|=L, \lb{5.53} \\
{\wti \cM}^\bot&=\{1,\dots,N\}\backslash\wti\cM \, \text{ for } \,
\wti \cM\subseteq \{1,\dots,N\}, \;  |\wti \cM|=M. \lb{5.54} 
\end{align}
Finally, if $\cL=\emptyset$ or $\cM=\emptyset$, then $K$ is a
Volterra operator and hence ${\det}_2  (I-K)=1$. 
\end{theorem}
\begin{proof}
Combining \eqref{3.29}, \eqref{5.42}, \eqref{5.43}, and \eqref{5.44} one
obtains
\begin{align}
{\det}_2 (I-K)&={\det}_{\bbC^M} (I_M-G)\exp\bigg(\int_0^\tau 
dx\, f_2(x)g_2(x)\bigg) \no\\
&={\det}_{\bbC^M} (I_M-G)\exp\bigg(\tau\sum_{m\in\cM}\beta_m\bigg) \no
\\
&={\det}_{\bbC^M} (I_M-G)\exp(\tau k(0_-)) \lb{5.55} \\
&=\exp\bigg(\tau k(0_-)-\tau\sum_{m\in\cM} \mu_m \bigg) 
\sum_{\substack{\wti \cL\subseteq \{1,\dots,N\}\\ |\wti \cL|=L}} 
V_{\wti\cL} \exp\bigg(-i\tau \sum_{m\in{\wti\cL}^\bot} \zeta_{m}\bigg), 
\no
\end{align}
where
\begin{align}
V_{\wti\cL}&=(-1)^M \bigg(\prod_{m\in{\wti\cL}^\bot}\beta_m\bigg) 
\bigg(\prod_{m'\in{\wti\cL}^\bot}\gamma_{m'}\bigg) 
\prod_{m''\in{\wti\cL}^\bot} \prod_{\substack{p\in{\wti\cL}^\bot
\\p\neq m''}} (i\zeta_{m''}-i\zeta_p) \lb{5.56} \\
& \quad \times\prod_{p'\in\cM} \prod_{\substack{p''\in\cM\\p''\neq p'}} 
(\mu_{p'}-\mu_{p''})\prod_{q\in{\wti\cL}^\bot}\prod_{q'\in\cM} 
(\mu_{q'}+i\zeta_q)^{-1}\prod_{r\in\cM}
\prod_{\substack{r'\in{\wti\cL}^\bot\\r'\neq r}}
(\mu_r+i\zeta_{r'})^{-1}. \no 
\end{align}
Elementary manipulations, using \eqref{5.20}, \eqref{5.21}, then reduce
\eqref{5.56} to \eqref{5.49} and hence prove \eqref{5.47}. To prove
\eqref{5.48} one can argue as follows. Introducing 
\begin{equation}
\wti{\cF(k)}(x)=\cF(k)(-x), \quad x\in\bbR \lb{5.57}
\end{equation}
with associated kernel function
\begin{equation}
\tilde k(t)=k(-t), \quad t\in\bbR\backslash\{0\}, \lb{5.58}
\end{equation}
equation \eqref{5.16} yields
\begin{equation}
1-\wti{\cF(k)}(x)=\prod_{n\in\cN} (x+\zeta_n)
\prod_{\ell\in\cL} (x-i\lambda_\ell)^{-1} 
\prod_{m\in\cM} (x+i\mu_m)^{-1}. \lb{5.59}
\end{equation}
Denoting by $\wti K$ the truncated Wiener--Hopf operator in
$L^2((0,\tau);dx)$ with convolution integral kernel $\tilde k$ (i.e., 
replacing $k$ by $\tilde k$ in \eqref{5.1}, and applying \eqref{5.47}
yields 
\begin{equation}
{\det}_2 (I-\wti K)=\exp\bigg(\tau \tilde k(0_-)-\tau\sum_{\ell\in\cL}
\lambda_\ell \bigg) \sum_{\substack{\wti \cM\subseteq 
\{1,\dots,N\}\\ |\wti \cM|=M}} W_{\wti\cM} \, \exp\big(i\tau
\sum_{\ell\in{\wti\cM}^\bot} \zeta_{\ell}\big). \lb{5.60} 
\end{equation}
Here $W_{\wti\cM}$ is given by \eqref{5.50} (after interchanging the 
roles of $\lambda_\ell$ and $\mu_m$ and interchanging $\zeta_m$ and
$-\zeta_\ell$, etc.) By \eqref{5.58}, $\tilde k(0_-)=k(0_+)$. Since
$\wti K=K^\top$, where $K^\top$ denotes the transpose integral operator
of $K$ (i.e., $K^\top$ has integral kernel $K(x',x)$ if $K(x,x')$ is the
integral kernel of $K$), and hence
\begin{equation}
{\det}_2(I-\wti K)={\det}_2(I-K^\top)={\det}_2(I- K), \lb{5.61}
\end{equation}
one arrives at \eqref{5.48}.

Finally, if $\cL=\emptyset$ then $k(0_+)=0$ and one infers 
${\det}_2 (I-K)=1$ by \eqref{5.48}. Similarly, if $\cM=\emptyset$, then
$k(0_-)=0$ and again ${\det}_2  (I-K)=1$ by \eqref{5.47}. 
\end{proof}

\begin{remark} \lb{r5.4}  
$(i)$ Theorem \ref{t5.3} permits some extensions. For instance, it
extends to the case where $\Re(\lambda_\ell) \geq 0$, $\Re(\mu_m)\geq 0$.
In this case the Fourier transform of $k$ should be understood in the
sense of distributions. One can also handle the case where
$-i\lambda_\ell$ and $i\mu_m$ are higher order poles of $\cF(k)$ by
using a limiting argument. \\
$(ii)$ The operator $K$ is a trace class operator, $K\in
\cB_1(L^2((0,\tau);dx))$, if and only if $k$ is continuous at $t=0$  
$($cf.\ equation $(2)$ on p.\ 267 and Theorem 10.3 in \cite{GGK96}$)$.  
\end{remark}

Explicit formulas for  determinants of Toeplitz operators with rational 
symbols are due to Day \cite{Da75}. Different proofs of Day's
formula can be found in \cite[Theorem 6.29]{BS83}, \cite{Go80}, and
\cite {HJ78}. Day's theorem requires that the degree of the
numerator of the rational symbol be greater or equal to that of the
denominator. An extension of Day's result avoiding such a restriction
recently appeared in \cite{CP01}.  Determinants of rationally generated
block operator matrices have also been studied in
\cite{Ti87} and \cite{Tr86}. Explicit representations for determinants
of the block-operator matrices of Toeplitz type with analytic symbol of
a special form has been obtained in \cite{Go92}. Textbook expositions of
these results can be found in \cite[Theorem 6.29]{BS83} and
\cite[Theorem 10.45]{BS90} (see also \cite[Sect.\ 5.9]{BS99}).

The explicit result \eqref{5.48}, that is, an explicit representation of
the $2$-modified Fredholm determinant for truncated Wiener-Hopf
operators on a finite interval, has first been obtained by B\"ottcher
\cite{Bo89}. He succceeded in reducing the problem to that of Toeplitz
operators combining a discretization approach and Day's formula. Theorem
\ref{t5.3} should thus be viewed as a continuous analog of Day's
formula. The method of proof presented in this paper based on
\eqref{3.29} is remarkably elementary and direct. A new method for the 
computation of ($2$-modified) determinants for truncated Wiener-Hopf
operators, based on the Nagy--Foias functional model, has recently been
suggested in {\cite{MP01} (cf.\ also {\cite{MP98}), without, however,
explicitly computing the right-hand sides of \eqref{5.47}, \eqref{5.48}.
A detailed exposition of the theory of operators of convolution type
with rational symbols on a finite interval, including representations
for resolvents, eigenfunctions, and (modified) Fredholm determinants
(different from the explicit one in Theorem \ref{t5.3}), can be found in
\cite[Sect.\ XIII.10]{GGK90}. Finally, extensions of the classical 
Szeg{\H o}--Kac--Achiezer formulas to the case of matrix-valued rational
symbols can be found in \cite{GK92} and \cite{GKS87}.

\medskip
\noindent {\bf Acknowledgements.} 
It is with great pleasure that we dedicate this paper to Eduard R.\
Tsekanovskii on the occasion of his 65th birthday. His contributions to
operator theory are profound and long lasting. In addition, we greatly
appreciate his beaming personality and, above all, his close friendship. 

We thank Radu Cascaval, David Cramer,
Vadim Kostrykin, Yuri Latushkin, and Barry Simon for useful discussions. 


\end{document}